\documentclass[3p,sort&compress,11pt,fleqn,review]{elsarticle}
\usepackage{setspace,tikz,float,lmodern,booktabs,amsmath,amssymb,siunitx,structmech,gnuplot-lua-tikz,subcaption}
\usetikzlibrary{shapes,arrows.meta}
\journal{}
\newcommand*{\figref}[1]{Fig.~\ref{#1}}
\newcommand*{\eqsref}[1]{Eq.~(\ref{#1})}

\newcommand*{\mb}[1]{\boldsymbol{#1}}
\newcommand*{\mT}{\mathrm{T}}
\newcommand*{\md}[1]{\mathrm{d}#1}
\newcommand*{\pfrac}[2]{\dfrac{\partial#1}{\partial#2}}
\tikzset{FFIXED/.style={postaction={draw,decorate,decoration={border,angle=-60,amplitude=2mm,segment length=1.9mm}}}}
\begin{document}
\begin{abstract}
This work presents a general finite element formulation based on a six--field variational principle that incorporates the consistent couple stress theory. A simple, efficient and local iteration free solving procedure that covers both elastic and inelastic materials is derived to minimise computation cost. With proper interpolations, membrane elements of various nodes are proposed as the examples. The implemented finite elements are used to conduct numerical experiments to investigate the performance of the in-plane drilling degrees of freedom introduced by the consistent couple stress theory. The mesh dependency issue is also studied with both elastic and inelastic materials. It is shown that the consistent couple stress theory provides an objective definition of rotation compared with the Cauchy theory but additional regularisation (or other techniques) is required to overcome mesh/size dependency in softening or fracture related problems. In the case of hardening continuum problems and/or large characteristic lengths, the proposed formulation and elements offer a more reliable approach to model structures with both translational and rotational degrees of freedom.
\end{abstract}
\begin{frontmatter}
\title{New mixed formulation and mesh dependency of finite elements based on the consistent couple stress theory}
\author[add1]{Theodore~L.~Chang\corref{tlc}}\ead{tlcfem@gmail.com}
\author[add2]{Chin-Long~Lee}
\cortext[tlc]{corresponding author}
\address[add1]{IRIS Adlershof, Humboldt-Universität zu Berlin, Berlin, Germany, 12489.}
\address[add2]{Department of Civil and Natural Resources Engineering, University of Canterbury, Christchurch, New Zealand, 8041.}
\end{frontmatter}
\section{Introduction}
The Cauchy continuum mechanics has received wide recognition over the years and has become the standard framework for many engineering disciplines. For simple, idealised problems, analytical solutions can be constructed \citep[see, e.g.,][]{Timoshenko2010} with proper assumptions. The Cauchy theory also provides a theoretical basis for the extensively used finite element methods. However, due to the lack of proper measure of rotation/curvature (and its force/stress conjugate), the Cauchy theory is incapable of describing the size effect, which is frequently seen in many materials at various scales (from nano to macro) and plays a vital role in fracture mechanics \citep{Bazant1984}. Besides, certain difficulties (e.g., kinematic compatibility among different types of elements) are encountered in problems involving elements with/without rotation field. A typical example would be the beam--panel joint in which the beam possesses in-plane rotation while the panel, often modelled by membrane elements, does not. To overcome the shortcomings of the Cauchy theory, researchers have been investigating alternatives, which are often called nonlocal theories. The justification of nonlocal continuum theories is also well explained by \citet{Eringen2004}. According to whom, several different formulations, for example, micropolar theory, micromorphic theory and gradient theory, are available for various applications. Recent advances and reviews on relevant topics can be seen elsewhere, \citep[to name a few,][]{Maugin2010,Altenbach2011,Srinivasa2017,Thai2017,Faghidian2018,Shaat2017,Aifantis2011,Salehipour2015,Apuzzo2018,Lim2015,Shaat2020}.

One popular branch, which belongs to the second gradient theory and utilises gradient of strain to define rotation, has received popularity over the recent years. A representative formulation proposed by \citet{Yang2002}, named as the modified couple stress theory, stands out and is widely used in various beam/plate related problems \citep[see, e.g.,][]{Khakalo2018,Farokhi2018,Kim2019,Thanh2019,Thanh2019a,Fan2020,Yuan2020,Lu2017}. However, as extensively compared and discussed by \citet{Hadjesfandiari2016}, the resulting formulation possesses some inconsistencies regarding boundary condition handling. Due to its inherit construction (coupling between in-plane and out-of-plane components), it is also difficult to be applied to plane stress problems. Among other proposals, although disputed by others \citep{Neff2016}, the consistent couple stress theory \citep{Hadjesfandiari2011}, an extension of the work by \citet{Mindlin1968}, provides a concise and promising framework by introducing curvature and couple stress vectors into the formulation. For further reviews on theoretical developments and experimental investigations regarding couple stress theories, interested readers can refer to the work by \citet{Pedgaonkar2021} and the references therein.

Due to the additional curvature introduced, the couple stress theory allows construction of both in-plane and out-of-plane rotation fields. There is no doubt that in-plane rotations can be constructed within the Cauchy framework. Early attempts date back to the work by \citet{Allman1984}. Recent research work on membrane elements with in-plane drilling degrees of freedom is reviewed by \citet{Boutagouga2020}. However, the drilling degree of freedom defined within the Cauchy theory shows significant mesh dependency, which limits its application in modelling problems involving, for example, both membrane and beam elements \citep{Chang2020}. This is, interestingly, not well discussed by researchers and drilling DoFs serve more like additional auxiliary DoFs to enrich existing displacement fields. In the couple stress theory, although the corresponding rigid body like rotation measures are defined as functions of translational displacements, by using either the Lagrange multiplier method or the penalty function method, it is possible to interpolate rotation fields independently. This provides additional degrees of freedom which can be used to construct, for example, shell and membrane elements with in-plane rotations.

Based on the above motivations, it would thus be of great interest to investigate how the drilling degrees of freedom defined in the consistent couple stress theory would perform subjected to moment/rotation like loads and whether the associated mesh size dependency could be alleviated. This paper aims to discuss this issue. To this end, a new general variational framework that is flexible enough to accommodate various interpolation schemes with the assist of the consistent couple stress theory \citep{Hadjesfandiari2011}. This paper is organised as follows. The consistent couple stress theory is briefly reviewed by summarising key equations and quantities that are later used in finite element formulation. Based on this theory, a general mixed variational principle involving six independent fields is derived. After formulating a system of linear equations, a simple, highly efficient, local iteration free solving strategy is designed for general inelastic material models. As the examples, two membrane elements, named as couple stress triangle (CSMT3) and quadrilateral (CSMQ4), are implemented by choosing proper interpolations. The performance of the proposed elements is further investigated by conducting numerical experiments with the focus on the performance regarding drilling degrees of freedom and its mesh dependency.
\section{The Consistent Couple Stress Theory}
Details of the adopted consistent couple stress theory can be seen elsewhere \citep{Hadjesfandiari2011}. Here a brief summary of all necessary equations later utilised in finite element formulation is presented. Compared with the original rigorous theory, here the definition of stress tensor is changed to allow a clearer formulation. As a result, the equilibrium takes a concise form. Readers shall be aware of the difference.
\subsection{Kinematics}
Within the infinitesimal strain framework, the couple stress theory \citep{Hadjesfandiari2011} accounts for both symmetric and skew--symmetric parts of displacement gradient in kinematics, that is
\begin{gather}\label{eq:varepsilon}
\varepsilon_{ij}=\dfrac{1}{2}\left(u_{i,j}+u_{j,i}\right),\\
\omega_{ij}=\dfrac{1}{2}\left(u_{i,j}-u_{j,i}\right),
\end{gather}
where $u_i$ is the displacement field, $\varepsilon_{ij}$ is the infinitesimal strain tensor, $\omega_{ij}$ is the skew--symmetric strain tensor that can be equivalently represented by the rotation vector $\theta_i$,
\begin{gather}\label{eq:theta}
\theta_i=\dfrac{1}{2}\epsilon_{ijk}\omega_{kj}=\dfrac{1}{2}\epsilon_{ijk}u_{k,j},
\end{gather}
where $\epsilon_{ijk}$ is the Levi--Civita (permutation) symbol. By utilising $\theta_i$, the mean curvature vector $\kappa_i$ is defined as
\begin{gather}\label{eq:kappa}
\kappa_i=\dfrac{1}{2}\epsilon_{ijk}\theta_{k,j},
\end{gather}
which can be further expressed in terms of $u_i$ so that
\begin{gather}
\kappa_i=\dfrac{1}{4}\left(u_{j,ij}-u_{i,jj}\right).
\end{gather}
The engineering mean curvature $k_i$ can be defined accordingly to be
\begin{gather}
k_i=-2\kappa_i.
\end{gather}
\subsection{Equilibrium}
The total stress tensor consists of two parts, namely the conventional symmetric stress tensor $\sigma_{ij}$ and the additional skew--symmetric couple stress tensor $\mu_{ij}$. Compared with the original literature \citep{Hadjesfandiari2011}, here a slightly different notation is used for brevity.

The couple stress vector $\mu_i$ dual to $\mu_{ij}$ can be defined as
\begin{gather}
\mu_{ji}=\epsilon_{ijk}\mu_k.
\end{gather}
It is found that $\mu_i$ is the energetic conjugate to engineering mean curvature $k_i$.

The corresponding equilibrium equation can be shown as
\begin{gather}\label{eq:equilibrium}
\sigma_{ji,j}+\dfrac{1}{2}\left(\mu_{j,ij}-\mu_{i,jj}\right)+f_i=0,
\end{gather}
where $f_i$ is the body force field.
\subsection{Isotropic Linear Elasticity}
The stored energy function $W$ can be defined as a function of $\varepsilon_{ij}$ and $\kappa_i$ such that
\begin{gather}\label{eq:stored_energy}
W\left(\varepsilon_{ij},\kappa_i\right)=\dfrac{1}{2}\lambda\varepsilon_{ii}\varepsilon_{jj}+\mu\varepsilon_{ij}\varepsilon_{ij}+8\eta\kappa_i\kappa_i,
\end{gather}
where $\lambda$ and $\mu$ are Lam\'{e} constants, $\eta$ is the additional material constant. The corresponding constitutive relations for linear elasticity can then be derived to be
\begin{gather}\label{eq:constitutive_couple}
\mu_i=-8\eta\kappa_i=4\eta{}k_i,\\
\sigma_{ij}=\lambda\delta_{ij}\varepsilon_{kk}+2\mu\varepsilon_{ij},
\end{gather}
where $\delta_{ij}$ is the Kronecker delta. The additional constant $\eta$ can be related to shear modulus $\mu$ by the characteristic length $l$ according to the following expression,
\begin{gather}
\eta=l^2\mu.
\end{gather}
It is worth mentioning that $l$ may not be a material constant and depends on the geometry of structures \citep{Khorshidi2018}. Its determination is thus not further extended in this work.

Given that in this work $\sigma_{ij}$ is assumed to be symmetric, without loss of generality, the variation of $W$ can be expressed as
\begin{gather}\label{eq:potential_energy}
\delta{}W\left(\varepsilon_{ij},\kappa_i\right)=\bar{\sigma}_{ij}\delta\varepsilon_{ij}-2\bar{\mu}_{i}\delta\kappa_i,
\end{gather}
or
\begin{gather}
\delta{}W\left(\varepsilon_{ij},k_i\right)=\bar{\sigma}_{ij}\delta\varepsilon_{ij}+\bar{\mu}_{i}\delta{}k_i,
\end{gather}
with $\bar{\sigma}_{ij}$ and $\bar{\mu}_{i}$ denoting stress and couple stress tensors obtained from typical strain driven constitutive model. For linear elasticity, they are simply
\begin{gather}
\bar{\sigma}_{ij}=\sigma_{ij},\qquad\bar{\mu}_{i}=\mu_{i}.
\end{gather}
\subsection{Remarks}
It could be noted that $C^1$ continuity is required by the displacement field $u_i$ due to the presence of the second order derivatives in $\kappa_i$. Shape functions based on functions, such as NURBS \citep[see][]{Dargush2021}, that support $C^1$ continuity can be adopted to construct proper finite elements. Alternatively, $\theta_i$ can be treated as an independent field and the corresponding kinematic equations can be introduced into the functional via the method of Lagrange multiplier, applications of which can be seen for example in the work by \citet{Darrall2013,Deng2016,Pedgaonkar2021}. Further discussions of such a consistent couple stress theory can also be seen in the work by \citet{Hadjesfandiari2016}.
\section{A Six--Field Mixed Framework}
In this section, the variational theorem developed by \citet{Darrall2013} is extended in this study to a more general form which resembles the Hu--Washizu variational theorem in the classic Cauchy theory. The approach follows a typical Lagrange multiplier method.

By assuming that the essential boundary conditions can be satisfied by proper construction, in absence of body force and surface force/moment traction, the total potential energy functional over an arbitrary domain $V$ can be simply expressed as
\begin{gather}
\varPi\left(\varepsilon_{ij},\kappa_i\right)=\int_VW~\md{V}
\end{gather}
where $W\left(\varepsilon_{ij},\kappa_i\right)$ is defined in \eqsref{eq:potential_energy} for linear elastic material. Since $\varepsilon_{ij}\left(u_i\right)$ and $\kappa_i\left(u_i\right)$ are functions of $u_i$, $u_i$ is the only independent field in the above functional.

Now consider the case in which $u_i$, $\theta_i$, $\varepsilon_{ij}$, and $\kappa_i$ are all treated as independent fields, the kinematic equations \eqsref{eq:varepsilon}, \eqsref{eq:kappa} and \eqsref{eq:theta} need to be satisfied in a weak form. By introducing three Lagrange multipliers $\alpha_{ij}$, $\beta_i$ and $\gamma_i$, those equations can be appended to the above functional so that
\begin{multline}
\varPi\left(u_i,\theta_i,\varepsilon_{ij},\kappa_i,\alpha_{ij},\beta_i,\gamma_i\right)=\int_VW~\md{V}
+\int_V\alpha_{ij}\left(\varepsilon_{ij}-\dfrac{1}{2}\left(u_{i,j}+u_{j,i}\right)\right)~\md{V}\\
+\int_V\beta_i\left(\kappa_i-\dfrac{1}{2}\epsilon_{ijk}\theta_{k,j}\right)~\md{V}
+\int_V\gamma_i\left(\theta_i-\dfrac{1}{2}\epsilon_{ijk}u_{k,j}\right)~\md{V}.
\end{multline}
Taking the first variation leads to
\begin{gather}\label{eq:variation}
\begin{split}
\delta\varPi&=
\int_V\bar{\sigma}_{ij}\delta\varepsilon_{ij}~\md{V}
-\int_V2\bar{\mu}_{i}\delta\kappa_i~\md{V}\\&
+\int_V\delta{}\alpha_{ij}\left(\varepsilon_{ij}-\dfrac{1}{2}\left(u_{i,j}+u_{j,i}\right)\right)~\md{V}
+\int_V\alpha_{ij}\left(\delta{}\varepsilon_{ij}-\dfrac{1}{2}\left(\delta{}u_{i,j}+\delta{}u_{j,i}\right)\right)~\md{V}\\&
+\int_V\delta{}\beta_i\left(\kappa_i-\dfrac{1}{2}\epsilon_{ijk}\theta_{k,j}\right)~\md{V}
+\int_V\beta_i\left(\delta{}\kappa_i-\dfrac{1}{2}\epsilon_{ijk}\delta{}\theta_{k,j}\right)~\md{V}\\&
+\int_V\delta{}\gamma_i\left(\theta_i-\dfrac{1}{2}\epsilon_{ijk}u_{k,j}\right)~\md{V}
+\int_V\gamma_i\left(\delta{}\theta_i-\dfrac{1}{2}\epsilon_{ijk}\delta{}u_{k,j}\right)~\md{V}.
\end{split}
\end{gather}
By performing integration by parts and applying the divergence theorem, one can find
\begin{gather}
-\dfrac{1}{2}\int_V\alpha_{ij}\left(\delta{}u_{i,j}+\delta{}u_{j,i}\right)~\md{V}=\int_V\dfrac{\alpha_{ij,j}+\alpha_{ji,j}}{2}\delta{}u_i~\md{V}-\int_S\dfrac{\alpha_{ij}+\alpha_{ji}}{2}n_j\delta{}u_i~\md{S},\\
-\dfrac{1}{2}\int_V\epsilon_{ijk}\beta_i\delta\theta_{k,j}~\md{V}=\dfrac{1}{2}\int_V\epsilon_{ijk}\beta_{i,j}\delta\theta_k~\md{V}-\dfrac{1}{2}\int_S\epsilon_{ijk}\beta_in_j\delta\theta_k~\md{S},\\
-\dfrac{1}{2}\int_V\epsilon_{ijk}\gamma_i\delta{}u_{k,j}~\md{V}=\dfrac{1}{2}\int_V\epsilon_{ijk}\gamma_{i,j}\delta{}u_k~\md{V}-\dfrac{1}{2}\int_S\epsilon_{ijk}\gamma_in_j\delta{}u_k~\md{S}.
\end{gather}
Here the boundary $S$ is not further refined for simplicity.

Inserting the above expressions back to \eqsref{eq:variation} gives
\begin{gather}
\begin{split}
\delta\varPi&=
\int_V\bar{\sigma}_{ij}\delta\varepsilon_{ij}~\md{V}+\int_V\alpha_{ij}\delta{}\varepsilon_{ij}~\md{V}
-\int_V2\bar{\mu}_{i}\delta\kappa_i~\md{V}+\int_V\beta_i\delta{}\kappa_i~\md{V}\\&
+\int_V\delta{}\alpha_{ij}\left(\varepsilon_{ij}-\dfrac{1}{2}\left(u_{i,j}+u_{j,i}\right)\right)~\md{V}
+\int_V\delta{}\beta_i\left(\kappa_i-\dfrac{1}{2}\epsilon_{ijk}\theta_{k,j}\right)~\md{V}\\&
+\int_V\delta{}\gamma_i\left(\theta_i-\dfrac{1}{2}\epsilon_{ijk}u_{k,j}\right)~\md{V}
+\int_V\dfrac{\alpha_{ij,j}+\alpha_{ji,j}}{2}\delta{}u_i~\md{V}
+\dfrac{1}{2}\int_V\epsilon_{ijk}\gamma_{i,j}\delta{}u_k~\md{V}\\&
+\int_V\gamma_i\delta{}\theta_i~\md{V}+\dfrac{1}{2}\int_V\epsilon_{ijk}\beta_{i,j}\delta\theta_k~\md{V}-\dfrac{1}{2}\int_S\epsilon_{ijk}\beta_in_j\delta\theta_k~\md{S}\\&
-\int_S\dfrac{\alpha_{ij}+\alpha_{ji}}{2}n_j\delta{}u_i~\md{S}-\dfrac{1}{2}\int_S\epsilon_{ijk}\gamma_in_j\delta{}u_k~\md{S}.
\end{split}
\end{gather}

Since the variations $\delta{}u_i$, $\delta{}\theta_i$, $\delta{}\varepsilon_{ij}$, $\delta{}\kappa_i$, $\delta{}\alpha_{ij}$, $\delta{}\beta_i$ and $\delta{}\gamma_i$ are arbitrary, the stationary condition requires the following equations involving Lagrange multipliers to hold.
\begin{gather}\label{eq:langrange}
\bar{\sigma}_{ij}+\alpha_{ij}=0,\quad
-2\bar{\mu}_i+\beta_i=0,\quad
\alpha_{kj,j}+\alpha_{jk,j}+\epsilon_{ijk}\gamma_{i,j}=0,\quad
\gamma_k+\dfrac{1}{2}\epsilon_{ijk}\beta_{i,j}=0.
\end{gather}
It can be identified that $\alpha_{ij}=-\sigma_{ij}$, $\beta_i=2\mu_i$ and $\gamma_i=\epsilon_{ijk}\mu_{k,j}$ with $\sigma_{ij}$ and $\mu_i$ be independent fields. The third equation in \eqsref{eq:langrange} is essentially
\begin{gather}
\begin{split}
0=\left(\alpha_{kj}+\alpha_{jk}\right)_{,j}+\epsilon_{ijk}\epsilon_{imn}\mu_{m,nj}=\left(\alpha_{ji}+\alpha_{ij}\right)_{,j}+\left(\mu_{j,i}-\mu_{i,j}\right)_{,j},
\end{split}
\end{gather}
which is the stress equilibrium \eqsref{eq:equilibrium} in absence of body force $f_i$.

The functional in its general form is then
\begin{multline}
\varPi\left(u_i,\theta_i,\varepsilon_{ij},\sigma_{ij},\kappa_i,\mu_i\right)=\int_VW~\md{V}
-\int_V\sigma_{ij}\left(\varepsilon_{ij}-\dfrac{1}{2}\left(u_{i,j}+u_{j,i}\right)\right)~\md{V}\\
+\int_V2\mu_i\left(\kappa_i-\dfrac{1}{2}\epsilon_{ijk}\theta_{k,j}\right)~\md{V}
+\int_V\epsilon_{imn}\mu_{n,m}\left(\theta_i-\dfrac{1}{2}\epsilon_{ijk}u_{k,j}\right)~\md{V}.
\end{multline}
It shall be noted that all boundary terms are omitted for brevity. In vector/matrix form, it can also be written as
\begin{multline}\label{eq:functional}
\varPi\left(\mb{u},\mb{\theta},\mb{\varepsilon},\mb{\sigma},\mb{\kappa},\mb{\mu}\right)=\int_VW~\md{V}
+\int_V\mb{\sigma}^\mT\left(\nabla^s\mb{u}-\mb{\varepsilon}\right)~\md{V}\\
+\int_V2\mb{\mu}^\mT\left(\mb{\kappa}-\dfrac{1}{2}\nabla\times\mb{\theta}\right)~\md{V}
+\int_V\left(\nabla\times\mb{\mu}\right)^\mT\left(\mb{\theta}-\dfrac{1}{2}\nabla\times\mb{u}\right)~\md{V}.
\end{multline}
The symbol $\nabla^s\mb{u}$ is used to denote the result of $\left(u_{i,j}+u_{j,i}\right)/2$ expressed in the Voigt form.

\eqsref{eq:functional} occupies a form similar to that of the functional used in the Hu--Washizu principle. Fields $\mb{u}$, $\mb{\theta}$ and $\mb{\mu}$ require $C^0$ continuity, while $\mb{\kappa}$, $\mb{\varepsilon}$ and $\mb{\sigma}$ can be constant fields. Starting from \eqsref{eq:functional}, various levels of simplifications can be conducted to derive both mixed--type and hybrid--type finite elements. For example, any of \eqsref{eq:varepsilon}, \eqsref{eq:kappa} and \eqsref{eq:theta} can be satisfied in strong forms thus the corresponding terms can be omitted from the functional. The mixed functional used by \citet{Darrall2013}, which is
\begin{gather}
\varPi\left(\mb{u},\mb{\theta},\mb{\mu}\right)=\int_VW~\md{V}
+\int_V\left(\nabla\times\mb{\mu}\right)^\mT\left(\mb{\theta}-\dfrac{1}{2}\nabla\times\mb{u}\right)~\md{V}+\varPi_{b.t.}.
\end{gather}
can be obtained by enforcing \eqsref{eq:varepsilon} and \eqsref{eq:kappa} in strong forms. Besides, by considering \eqsref{eq:varepsilon} only, another functional can be obtained.
\begin{multline}
\varPi\left(\mb{u},\mb{\theta},\mb{\kappa},\mb{\mu}\right)=\int_VW~\md{V}
+\int_V2\mb{\mu}^\mT\left(\mb{\kappa}-\dfrac{1}{2}\nabla\times\mb{\theta}\right)~\md{V}\\
+\int_V\left(\nabla\times\mb{\mu}\right)^\mT\left(\mb{\theta}-\dfrac{1}{2}\nabla\times\mb{u}\right)~\md{V}+\varPi_{b.t.}.
\end{multline}
One can also apply the divergence theorem to terms involving derivatives to convert between volume and surface integrals.

If, instead of $\mb{\kappa}$, engineering mean curvature $\mb{k}$ is used, an equivalent form can be derived following the same procedure.
\begin{multline}\label{eq:functional2}
\varPi\left(\mb{u},\mb{\theta},\mb{\varepsilon},\mb{\sigma},\mb{k},\mb{\mu}\right)=\int_VW~\md{V}
+\int_V\mb{\sigma}^\mT\left(\nabla^s\mb{u}-\mb{\varepsilon}\right)~\md{V}\\
-\int_V\mb{\mu}^\mT\left(\mb{k}+\nabla\times\mb{\theta}\right)~\md{V}
+\int_V\left(\nabla\times\mb{\mu}\right)^\mT\left(\mb{\theta}-\dfrac{1}{2}\nabla\times\mb{u}\right)~\md{V}.
\end{multline}
The difference between \eqsref{eq:functional} and \eqsref{eq:functional2} is solely the term $\displaystyle\int_V2\mb{\mu}^\mT\mb{\kappa}~\md{V}=-\int_V\mb{\mu}^\mT\mb{k}~\md{V}$.
\section{Finite Element Formulation}
In this section, the linear equation system of the aforementioned general six--field variational principle \eqsref{eq:functional} is derived. Since there is no other local residual apart from the one due to potential inelastic constitutive models, a locally iterative algorithm is not required. The final elemental stiffness may possess a form similar to that of conventional displacement based elements.
\subsection{Linear System}
Let six fields be discretized as follows.
\begin{gather}
\mb{u}=\mb{\phi}_{\mb{u}}\mb{p},\quad
\mb{\theta}=\mb{\phi}_{\mb{\theta}}\mb{q},\quad
\mb{\kappa}=\mb{\phi}_{\mb{\kappa}}\mb{r},\quad
\mb{\mu}=\mb{\phi}_{\mb{\mu}}\mb{s},\quad
\mb{\varepsilon}=\mb{\phi}_{\mb{\varepsilon}}\mb{\beta},\quad
\mb{\sigma}=\mb{\phi}_{\mb{\sigma}}\mb{\alpha}.
\end{gather}
Then, naturally, $\nabla^s\mb{u}=\mb{L}\mb{\phi}_{\mb{u}}\mb{p}$ where $\mb{L}$ is the differential operator which can be expressed as
\begin{gather}
\mb{L}=\begin{bmatrix}
\pfrac{}{x}&\cdot&\cdot&\pfrac{}{y}&\cdot&\pfrac{}{z}\\[4mm]
\cdot&\pfrac{}{y}&\cdot&\pfrac{}{x}&\pfrac{}{z}&\cdot\\[4mm]
\cdot&\cdot&\pfrac{}{z}&\cdot&\pfrac{}{y}&\pfrac{}{x}
\end{bmatrix}^\mT
\end{gather}
in 3D space. Similarly, the curl operator can be expressed as $\dfrac{1}{2}\nabla\times\left(\cdot\right)=\mb{J}\left(\cdot\right)$ with
\begin{gather}
\mb{J}=\dfrac{1}{2}\begin{bmatrix}
\cdot&-\pfrac{}{z}&\pfrac{}{y}\\[4mm]
\pfrac{}{z}&\cdot&-\pfrac{}{x}\\[4mm]
-\pfrac{}{y}&\pfrac{}{x}&\cdot
\end{bmatrix}.
\end{gather}

Now the functional can be rewritten as
\begin{multline}\label{eq:functional_new}
\varPi\left(\mb{p},\mb{q},\mb{r},\mb{s},\mb{\beta},\mb{\alpha}\right)=\int_VW~\md{V}
+\int_V\mb{\alpha}^\mT\mb{\phi}_{\mb{\sigma}}^\mT\left(\mb{L}\mb{\phi}_{\mb{u}}\mb{p}-\mb{\phi}_{\mb{\varepsilon}}\mb{\beta}\right)~\md{V}\\
+\int_V2\mb{s}^\mT\mb{\phi}_{\mb{\mu}}^\mT\left(\mb{\phi}_{\mb{\kappa}}\mb{r}-\mb{J}\mb{\phi}_{\mb{\theta}}\mb{q}\right)~\md{V}
+\int_V2\mb{s}^\mT\left(\mb{J}\mb{\phi}_{\mb{\mu}}\right)^\mT\left(\mb{\phi}_{\mb{\theta}}\mb{q}-\mb{J}\mb{\phi}_{\mb{u}}\mb{p}\right)~\md{V}.
\end{multline}

Taking variations leads to the following system of equations.
\begin{gather}\label{eq:stationary}
\left\{
\begin{array}{lll}
\dfrac{\delta\varPi}{\delta{}\mb{p}}=\mb{0}&\longrightarrow&\displaystyle\int_V\left(\mb{L}\mb{\phi}_{\mb{u}}\right)^\mT\mb{\phi}_{\mb{\sigma}}\mb{\alpha}-2\left(\mb{J}\mb{\phi}_{\mb{u}}\right)^\mT\left(\mb{J}\mb{\phi}_{\mb{\mu}}\right)\mb{s}~\md{V}=\mb{R}_{\mb{u}},\\[4mm]
\dfrac{\delta\varPi}{\delta{}\mb{q}}=\mb{0}&\longrightarrow&\displaystyle\int_V2\left(\mb{\phi}_{\mb{\theta}}^\mT\mb{J}\mb{\phi}_{\mb{\mu}}-\left(\mb{J}\mb{\phi}_{\mb{\theta}}\right)^\mT\mb{\phi}_{\mb{\mu}}\right)\mb{s}~\md{V}=\mb{R}_{\mb{\theta}},\\[4mm]
\dfrac{\delta\varPi}{\delta{}\mb{r}}=\mb{0}&\longrightarrow&\displaystyle\int_V\mb{\phi}_{\mb{\kappa}}^\mT{}W_{\mb{\kappa}}+2\mb{\phi}_{\mb{\kappa}}^\mT\mb{\phi}_{\mb{\mu}}\mb{s}~\md{V}=\mb{0},\\[4mm]
\dfrac{\delta\varPi}{\delta{}\mb{s}}=\mb{0}&\longrightarrow&\displaystyle
\int_V2\left(\left(\mb{J}\mb{\phi}_{\mb{\mu}}\right)^\mT\mb{\phi}_{\mb{\theta}}-\mb{\phi}_{\mb{\mu}}^\mT\mb{J}\mb{\phi}_{\mb{\theta}}\right)\mb{q}+2\mb{\phi}_{\mb{\mu}}^\mT\mb{\phi}_{\mb{\kappa}}\mb{r}-2\left(\mb{J}\mb{\phi}_{\mb{\mu}}\right)^\mT\mb{J}\mb{\phi}_{\mb{u}}\mb{p}~\md{V}=\mb{0},\\[4mm]
\dfrac{\delta\varPi}{\delta{}\mb{\beta}}=\mb{0}&\longrightarrow&\displaystyle\int_V\mb{\phi}_{\mb{\varepsilon}}^\mT{}W_{\mb{\varepsilon}}~\md{V}
-\int_V\mb{\phi}_{\mb{\varepsilon}}^\mT\mb{\phi}_{\mb{\sigma}}\mb{\alpha}~\md{V}=\mb{0},\\[4mm]
\dfrac{\delta\varPi}{\delta{}\mb{\alpha}}=\mb{0}&\longrightarrow&\displaystyle\int_V\mb{\phi}_{\mb{\sigma}}^\mT\mb{L}\mb{\phi}_{\mb{u}}\mb{p}-\mb{\phi}_{\mb{\sigma}}^\mT\mb{\phi}_{\mb{\varepsilon}}\mb{\beta}~\md{V}=\mb{0}.
\end{array}
\right.
\end{gather}
In the above system, $\mb{R}_{\mb{u}}$ and $\mb{R}_{\mb{\theta}}$ are nodal forces/resistances due to omitted boundary terms, $W_{\mb{\kappa}}$ and $W_{\mb{\varepsilon}}$ denote the partial derivatives respectively.

By further denoting
\begin{gather}
\mb{E}_1=\int_V\mb{\phi}_{\mb{\varepsilon}}^\mT\mb{C}\mb{\phi}_{\mb{\varepsilon}}~\md{V},\quad
\mb{E}_2=\int_V\mb{\phi}_{\mb{\kappa}}^\mT\mb{D}\mb{\phi}_{\mb{\kappa}}~\md{V},\quad
\mb{H}_1=-2\int_V\left(\mb{J}\mb{\phi}_{\mb{u}}\right)^\mT\left(\mb{J}\mb{\phi}_{\mb{\mu}}\right)~\md{V},\\
\mb{H}_2=\int_V\left(\mb{L}\mb{\phi}_{\mb{u}}\right)^\mT\mb{\phi}_{\mb{\sigma}}~\md{V},\quad
\mb{H}_3=2\int_V\mb{\phi}_{\mb{\theta}}^\mT\mb{J}\mb{\phi}_{\mb{\mu}}-\left(\mb{J}\mb{\phi}_{\mb{\theta}}\right)^\mT\mb{\phi}_{\mb{\mu}}~\md{V},\\
\mb{H}_4=-2\int_V\mb{\phi}_{\mb{\kappa}}^\mT\mb{\phi}_{\mb{\mu}}~\md{V},\quad
\mb{H}_5=\int_V\mb{\phi}_{\mb{\varepsilon}}^\mT\mb{\phi}_{\mb{\sigma}}~\md{V},\quad
\end{gather}
in which $\mb{C}$ and $\mb{D}$ denote material tangent moduli, the incremental form of linear system can be expressed as
\begin{gather}\label{eq:incremental_form}
\begin{bmatrix}
\cdot&\cdot&\cdot&\mb{H}_1&\cdot&\mb{H}_2\\
\cdot&\cdot&\cdot&\mb{H}_3&\cdot&\cdot\\
\cdot&\cdot&\mb{E}_2&-\mb{H}_4&\cdot&\cdot\\
\mb{H}_1^\mT&\mb{H}_3^\mT&-\mb{H}_4^\mT&\cdot&\cdot&\cdot\\
\cdot&\cdot&\cdot&\cdot&\mb{E}_1&-\mb{H}_5\\
\mb{H}_2^\mT&\cdot&\cdot&\cdot&-\mb{H}_5^\mT&\cdot
\end{bmatrix}
\begin{bmatrix}
\Delta\mb{p}\\\Delta\mb{q}\\\Delta\mb{r}\\\Delta\mb{s}\\\Delta\mb{\beta}\\\Delta\mb{\alpha}
\end{bmatrix}=\mb{R}.
\end{gather}
In which, $\mb{R}$ is used to denote residual whose discrete form is not shown. As can be seen later, it is possible to compute resistance directly based on \eqsref{eq:stationary}.
\subsection{Solution Procedure}
\subsubsection{Option One --- Invertible $\mb{E}_1$ and $\mb{E}_2$}
Since $\mb{E}_1$ and $\mb{E}_2$ are guaranteed to be square matrices, typically they are invertible with sufficient integration points provided, performing static condensation leads to the following elemental stiffness.
\begin{gather}
\mb{K}=
\begin{bmatrix}
\mb{H}_2\left(\mb{H}_5^\mT\mb{E}_1^{-1}\mb{H}_5\right)^{-1}\mb{H}_2^\mT+\mb{H}_1\left(\mb{H}_4^\mT\mb{E}_2^{-1}\mb{H}_4\right)^{-1}\mb{H}_1^\mT&\mb{H}_1\left(\mb{H}_4^\mT\mb{E}_2^{-1}\mb{H}_4\right)^{-1}\mb{H}_3^\mT\\
\mb{H}_3\left(\mb{H}_4^\mT\mb{E}_2^{-1}\mb{H}_4\right)^{-1}\mb{H}_1^\mT&\mb{H}_3\left(\mb{H}_4^\mT\mb{E}_2^{-1}\mb{H}_4\right)^{-1}\mb{H}_3^\mT
\end{bmatrix}.
\end{gather}
The above procedure requires matrix inversions to be conducted at element level whenever tangent moduli update. This is considered not very efficient for nonlinear analysis.
\paragraph{Stability Condition}
Let the size of $\mb{\phi}_{\mb{\varepsilon}}$ be $i\times{}j$ and the size of $\mb{\phi}_{\mb{\sigma}}$ be $i\times{}k$, then the number of integration points $n_{ip}$ must satisfy the following inequality for $\mb{E}_1$ to be invertible,
\begin{gather}
n_{ip}\cdot{}\min\left(i,~j\right)\geqslant{}\max\left(i,~j\right).
\end{gather}
Furthermore, for $\mb{H}_5^\mT\mb{E}_1^{-1}\mb{H}_5$ to be invertible,
\begin{gather}
n_{ip}\cdot{}\min\left(i,~j,~k\right)\geqslant{}\max\left(i,~j,~k\right).
\end{gather}
This, however, only guarantees successful computation of the elemental stiffness matrix thus is the minimum condition. Denoting the number of rigid body motion as $n_f$ and the size of elemental stiffness matrix as $n_k$, then for a single element, at least
\begin{gather}
n_{ip}\cdot{}\min\left(i,~j,~k\right)\geqslant{}n_k-n_f.
\end{gather}
Similar procedure can be applied to term $\mb{H}_4^\mT\mb{E}_2^{-1}\mb{H}_4$. Let the size of $\mb{\phi}_{\mb{\kappa}}$ be $l\times{}m$ and the size of $\mb{\phi}_{\mb{\mu}}$ be $l\times{}n$, then
\begin{gather}
n_{ip}\cdot{}\min\left(l,~m,~n\right)\geqslant{}\max\left(l,~m,~n,~n_k-n_f\right).
\end{gather}
Eventually,
\begin{gather}\label{eq:stability}
n_{ip}\cdot{}\min\left(i,~j,~k,~l,~m,~n\right)\geqslant{}\max\left(i,~j,~k,~l,~m,~n,~n_k-n_f\right).
\end{gather}
\subsubsection{Option Two --- Invertible $\mb{H}_4$ and $\mb{H}_5$}
To improve computational efficiency, an alternative strategy is proposed as follows. By proper construction, for example, providing sufficient number of integration points and choosing the same shape functions for $\mb{\kappa}$ ($\mb{\varepsilon}$) and $\mb{\mu}$ ($\mb{\sigma}$), $\mb{H}_4$ and $\mb{H}_5$ can be square and invertible. Similar to the work by the authors \citep{Chang2019}, by static condensation, \eqsref{eq:incremental_form} can be transformed to the one,
\begin{gather}
\begin{bmatrix}
\mb{H}_2\mb{H}_5^{-1}\mb{E}_1\mb{H}_5^{-\mT}\mb{H}_2^\mT+\mb{H}_1\mb{H}_4^{-1}\mb{E}_2\mb{H}_4^{-\mT}\mb{H}_1^\mT&\mb{H}_1\mb{H}_4^{-1}\mb{E}_2\mb{H}_4^{-\mT}\mb{H}_3^\mT\\
\mb{H}_3\mb{H}_4^{-1}\mb{E}_2\mb{H}_4^{-\mT}\mb{H}_1^\mT&\mb{H}_3\mb{H}_4^{-1}\mb{E}_2\mb{H}_4^{-\mT}\mb{H}_3^\mT
\end{bmatrix}
\begin{bmatrix}
\Delta\mb{p}\\\Delta\mb{q}
\end{bmatrix}=\mb{R}.
\end{gather}

The above system can be conveniently converted into a more expressive form, in which the elemental stiffness matrix $\mb{K}$ can be shown as
\begin{gather}\label{eq:stiffness}
\mb{K}=\mb{K}_1+\mb{K}_2
\end{gather}
with
\begin{gather}
\mb{K}_1=\begin{bmatrix}
\mb{H}_2\mb{H}_5^{-1}&\cdot\\\cdot&\cdot
\end{bmatrix}\begin{bmatrix}
\mb{E}_1&\cdot\\\cdot&\cdot
\end{bmatrix}\begin{bmatrix}
\mb{H}_5^{-\mT}\mb{H}_2^\mT&\cdot\\\cdot&\cdot
\end{bmatrix},\\
\mb{K}_2=\begin{bmatrix}
\cdot&\mb{H}_1\mb{H}_4^{-1}\\\cdot&\mb{H}_3\mb{H}_4^{-1}
\end{bmatrix}\begin{bmatrix}
\cdot&\cdot\\\cdot&\mb{E}_2
\end{bmatrix}\begin{bmatrix}
\cdot&\cdot\\\mb{H}_4^{-\mT}\mb{H}_1^\mT&\mb{H}_4^{-\mT}\mb{H}_3^\mT
\end{bmatrix}.
\end{gather}
The first term $\mb{K}_1$ solely depends on conventional displacement, strain and stress as defined in the Cauchy theory. Thus it can be derived from the Hu--Washizu theorem.

Since all $\mb{H}_n$ matrices are constant once interpolations are determined, they only need to be computed once during the initialisation stage. Furthermore, $\mb{H}_n$ require no additional storage as $\mb{H}_2\mb{H}_5^{-1}\mb{\phi}_{\mb{\varepsilon}}^\mT$, $\mb{H}_1\mb{H}_4^{-1}\mb{\phi}_{\mb{\kappa}}^\mT$ and $\mb{H}_3\mb{H}_4^{-1}\mb{\phi}_{\mb{\kappa}}^\mT$ can be stored as `equivalent strain matrices' for each integration point. Once the elemental stiffness matrix is computed, reordering of degrees of freedom may be performed.

It can be noted that the elemental stiffness possesses a symmetric structure although moduli $\mb{E}_1$ and $\mb{E}_2$ may be asymmetric due to, for example, non-associative plasticity in the case of material nonlinearity.
\paragraph{Stability Condition}
This option only requires the inversion of square matrices $\mb{H}_4$ and $\mb{H}_5$. Thus $\mb{\phi}_{\mb{\varepsilon}}$ ($\mb{\phi}_{\mb{\kappa}}$) and $\mb{\phi}_{\mb{\sigma}}$ ($\mb{\phi}_{\mb{\mu}}$) have the same size. Let $i\times{}j$ be the size of $\mb{\phi}_{\mb{\varepsilon}}$ and $m\times{}n$ be the size of $\mb{\phi}_{\mb{\kappa}}$, then the stability condition becomes
\begin{gather}\label{eq:stability2}
n_{ip}\cdot{}\min\left(i,~j,~m,~n\right)\geqslant{}\max\left(i,~j,~m,~n,~n_k-n_f\right).
\end{gather}
\subsubsection{Resistance}
From the third equation of the stationary condition \eqsref{eq:stationary}, $\mb{s}$ solely depends on material response thus can be computed directly as
\begin{gather}\label{eq:s}
\mb{s}=\mb{H}_4^{-1}\int_V\mb{\phi}_{\mb{\kappa}}^\mT{}\bar{\mb{\mu}}~\md{V}.
\end{gather}
Applying the same strategy to the fifth equation, one can obtain
\begin{gather}\label{eq:alpha}
\mb{\alpha}=\mb{H}_5^{-1}\int_V\mb{\phi}_{\mb{\varepsilon}}^\mT{}\bar{\mb{\sigma}}~\md{V}.
\end{gather}
In the above equations, similar to \eqsref{eq:potential_energy}, $\bar{\mb{\sigma}}$ and $\bar{\mb{\mu}}$ denote the stress and couple stress obtained from material model. Inserting \eqsref{eq:s} and \eqsref{eq:alpha} back to the first two equations in \eqsref{eq:stationary}, one can obtain the elemental resistance as
\begin{gather}
\begin{bmatrix}
\mb{R}_{\mb{u}}\\\mb{R}_{\mb{\theta}}
\end{bmatrix}=\begin{bmatrix}
\displaystyle\mb{H}_2\mb{H}_4^{-1}\int_V\mb{\phi}_{\mb{\kappa}}^\mT{}\bar{\mb{\mu}}~\md{V}+\mb{H}_1\mb{H}_5^{-1}\int_V\mb{\phi}_{\mb{\varepsilon}}^\mT{}\bar{\mb{\sigma}}~\md{V}\\
\displaystyle\mb{H}_3\mb{H}_4^{-1}\int_V\mb{\phi}_{\mb{\kappa}}^\mT{}\bar{\mb{\mu}}~\md{V}
\end{bmatrix}.
\end{gather}
Once material response is computed, elemental resistance can be conveniently assembled without worrying about local residual in each iteration.
\subsubsection{Remarks}
Two solution procedures are presented and compared. The first option, as the conventional procedure, requires heavy matrix inversions to compute the elemental stiffness matrix while the proposed alternative does not. Thus the second option possesses less space and time complexity and thus is considerably more computationally efficient.

Nevertheless, both options can be chosen and they share similar stability conditions. However, they may require different numbers of shape functions for the corresponding fields in order to construct the square matrices. Depending on the type of element and the number of available integration points, one may switch from one option to another for optimal formulation (mostly like the proposed one).

If \eqsref{eq:functional2} is used, except for replacing all $\mb{\phi}_{\mb{\kappa}}$ with $\mb{\phi}_{\mb{k}}$, only $\mb{H}_4$ needs to be modified.
\begin{gather}
\mb{H}_4=\int_V\mb{\phi}_{\mb{k}}^\mT\mb{\phi}_{\mb{\mu}}~\md{V},
\end{gather}
where $\mb{\phi}_{\mb{k}}$ is the interpolation of engineering mean curvature $\mb{k}$.
\subsection{Some Basic Elements}
\subsubsection{An Elementary Triangular Membrane Element --- CSMT3}
Not all couple stress theories support membrane problems, the popular modified couple stress theory \citep{Yang2002} adopts a symmetric couple stress tensor, with which in-plane response cannot be fully decoupled from out-of-plane response. Discussions of relevant topics can be seen elsewhere \citep{Hadjesfandiari2016}. The consistent couple stress theory is free from similar issues, which makes it more appealing for a wide range of general continuum problems.

The simplest membrane element may be the three-node triangular element. For plane stress problems, given the constitutive equation between mean curvature vector $\kappa_i$ and couple stress vector $\mu_i$ occupies the form shown in \eqsref{eq:constitutive_couple}, it is clear that in-plane and out-of-plane actions are decoupled. Thus six fields reduce to the following Voigt forms.
\begin{gather}
\mb{u}=\begin{bmatrix}
u_x\\u_y
\end{bmatrix},\quad
\mb{\theta}=\begin{bmatrix}
\theta_z
\end{bmatrix},\quad
\mb{\kappa}=\begin{bmatrix}
\kappa_x\\\kappa_y
\end{bmatrix},\quad
\mb{\mu}=\begin{bmatrix}
\mu_x\\\mu_y
\end{bmatrix},\quad
\mb{\varepsilon}=\begin{bmatrix}
\varepsilon_x\\\varepsilon_y\\\gamma_{xy}
\end{bmatrix},\quad
\mb{\sigma}=\begin{bmatrix}
\sigma_x\\\sigma_y\\\tau_{xy}
\end{bmatrix}.
\end{gather}

The linear mapping is used for coordinate, displacement, drilling rotation and couple stress, that is
\begin{gather}
\chi=\sum_{i=1}^3N_i\chi_i
\end{gather}
where $N_i$ is the complete first order shape function with modes $\begin{bmatrix}
1&x&y
\end{bmatrix}$, $\chi$ represents any of $x$, $y$, $u_x$, $u_y$, $\theta_z$, $\kappa_x$, $\kappa_y$, $\mu_x$ and $\mu_y$. Thus for each field, three nodal values are used for interpolation.

Given that $\mb{u}$ is linearly interpolated, $\mb{\varepsilon}$ and $\mb{\sigma}$ can be chosen to be constant fields.
\begin{gather}
\mb{\varepsilon}=\begin{bmatrix}
1&\cdot&\cdot\\\cdot&1&\cdot\\\cdot&\cdot&1
\end{bmatrix}\begin{bmatrix}
\beta_1\\\beta_2\\\beta_3
\end{bmatrix},\qquad
\mb{\sigma}=\begin{bmatrix}
1&\cdot&\cdot\\\cdot&1&\cdot\\\cdot&\cdot&1
\end{bmatrix}\begin{bmatrix}
\alpha_1\\\alpha_2\\\alpha_3
\end{bmatrix}.
\end{gather}
In this case, $\beta_i$ and $\alpha_i$ are essentially strain and stress components.

The explicit forms are not listed here for brevity. Interested readers are referred to typical textbooks on finite element methods for details of formulating interpolation matrices \citep[see, e.g.,][section 5.1.3.1]{Zienkiewicz2013}. Three integration points are used for numerical integration. This element is denoted as the couple stress mixed triangle (CSMT3) element.

The CSMT3 element is essentially an extension of the constant strain triangle (CST) element (CPS3 and CPE3 in ABAQUS notion) with additional fields incorporating mean curvature and couple stress. Given that only linear and constant interpolations are used, its performance shall be similar to that of the mixed version of CST.
\subsubsection{A Four-Node Quadrilateral Membrane Element --- CSMQ4}
Similarly, the aforementioned formulation can be applied to four-node quadrilaterals. The standard isoparametric mapping can be applied.
\begin{gather}
\chi=\sum_{i=1}^4N_i\chi_i
\end{gather}
with $N_i=\left(1+\xi_i\xi\right)\left(1+\eta_i\eta\right)/4$ where $\xi$ and $\eta$ are parent coordinates while $\xi_i$ and $\eta_i$ are values of parent coordinates of target node. Again, $\chi$ can be any of $x$, $y$, $u_x$, $u_y$, $\theta_z$, $\kappa_x$, $\kappa_y$, $\mu_x$ and $\mu_y$. A complete first order interpolation can be chosen for both strain and stress.
\begin{gather}
\mb{\phi}_{\mb{\varepsilon}}=\mb{\phi}_{\mb{\sigma}}=\begin{bmatrix}
1&\cdot&\cdot&x&\cdot&\cdot&y&\cdot&\cdot\\
\cdot&1&\cdot&\cdot&x&\cdot&\cdot&y&\cdot\\
\cdot&\cdot&1&\cdot&\cdot&x&\cdot&\cdot&y
\end{bmatrix}.
\end{gather}
It shall be mentioned that in the consistent couple stress theory, $\mb{\sigma}$ is not only governed by the Airy stress function but also an additional stress function \citep{Hadjesfandiari2011} from which $\mb{\mu}$ can be derived. Based on the previous discovery \citep{Chang2019}, $\mb{\phi}_{\mb{\varepsilon}}$ can be further modified to include the Poisson effect, which overcomes (near) incompressible problems. A \num{2x2} Gauss quadrature is used for numerical integration, this element is denoted as the couple stress mixed quadrilateral (CSMQ4) element.

Given that the displacement interpolation only uses four modes $\begin{bmatrix}
1&x&y&xy
\end{bmatrix}$, one can further define internal degrees of freedom to complete the second order interpolation so that performance can be potentially improved.
\subsection{Other Elements}
Since the incremental form and elemental stiffness are given in \eqsref{eq:incremental_form} and \eqsref{eq:stiffness}, there is no difficulty in constructing various types of elements such as serendipity quadrilaterals with various nodes, axisymmetric triangle/quadrilateral, eight--node brick/cube, four--node tetrahedron and higher order elements. The procedure follows a standard FEM approach. Apart from the above two elements (CSMT3 and CSMQ4), other elements implemented \citep{Chang2021} include: six--node second order triangle (CSMT6), serendipity quadrilaterals with five to eight nodes (CSMQ5, CSMQ6, CSMQ7, CSMQ8). The following table summarises the specifications of those elements. The quadrature schemes listed are the ones used in implementation, they can be replaced by other schemes depending on different interpolations.
\begin{table}[H]
\centering\footnotesize\caption{summary of implemented elements}\label{tab:summary}
\begin{tabular}{m{2cm}m{2cm}m{6cm}m{3cm}}
	\toprule
	designation & quadrature             & description                           & illustration      \\ \midrule
	CSMT3       & three-point            & three--node first order triangle      & \begin{tikzpicture}
\clip(-1,-1)rectangle(1,1);
\draw[line width=.4mm](-.8,-.7)node[circle,fill=black]{}--(.8,-.2)node[circle,fill=black]{}--(-.1,.8)node[circle,fill=black]{}--cycle;
\end{tikzpicture} \\
	CSMT6       & four-point             & six--node second order triangle       & \begin{tikzpicture}
\clip(-1,-1)rectangle(1,1);
\draw[line width=.4mm](-.8,-.7)node[circle,fill=black]{}--(.8,-.2)node[circle,fill=black,midway]{}node[circle,fill=black]{}--(-.1,.8)node[circle,fill=black,midway]{}node[circle,fill=black]{}--cycle node[circle,fill=black,midway]{};
\end{tikzpicture} \\
	CSMQ4       & \num{2x2} Gauss & four--node bilinear quadrilateral     & \begin{tikzpicture}
\clip(-1,-1)rectangle(1,1);
\draw[line width=.4mm](-.8,-.7)node[circle,fill=black]{}--(.8,-.5)node[circle,fill=black]{}--(.7,.6)node[circle,fill=black]{}--(-.7,.8)node[circle,fill=black]{}--cycle;
\end{tikzpicture} \\
	CSMQ5       & \num{3x3} Gauss & five--node serendipity quadrilateral  & \begin{tikzpicture}
\clip(-1,-1)rectangle(1,1);
\draw[line width=.4mm](-.8,-.7)node[circle,fill=black]{}--(.8,-.5)node[circle,fill=black,midway]{}node[circle,fill=black]{}--(.7,.6)node[circle,fill=black]{}--(-.7,.8)node[circle,fill=black]{}--cycle;
\end{tikzpicture} \\
	CSMQ6       & \num{3x3} Gauss & six--node serendipity quadrilateral   & \begin{tikzpicture}
\clip(-1,-1)rectangle(1,1);
\draw[line width=.4mm](-.8,-.7)node[circle,fill=black]{}--(.8,-.5)node[circle,fill=black,midway]{}node[circle,fill=black]{}--(.7,.6)node[circle,fill=black]{}--(-.7,.8)node[circle,fill=black,midway]{}node[circle,fill=black]{}--cycle;
\end{tikzpicture} \\
	CSMQ7       & \num{3x3} Gauss & seven--node serendipity quadrilateral & \begin{tikzpicture}
\clip(-1,-1)rectangle(1,1);
\draw[line width=.4mm](-.8,-.7)node[circle,fill=black]{}--(.8,-.5)node[circle,fill=black,midway]{}node[circle,fill=black]{}--(.7,.6)node[circle,fill=black,midway]{}node[circle,fill=black]{}--(-.7,.8)node[circle,fill=black,midway]{}node[circle,fill=black]{}--cycle;
\end{tikzpicture} \\
	CSMQ8       & \num{3x3} Gauss & eight--node biquadratic quadrilateral & \begin{tikzpicture}
\clip(-1,-1)rectangle(1,1);
\draw[line width=.4mm](-.8,-.7)node[circle,fill=black]{}--(.8,-.5)node[circle,fill=black,midway]{}node[circle,fill=black]{}--(.7,.6)node[circle,fill=black,midway]{}node[circle,fill=black]{}--(-.7,.8)node[circle,fill=black,midway]{}node[circle,fill=black]{}--cycle node[circle,fill=black,midway]{};
\end{tikzpicture} \\ \bottomrule
\end{tabular}
\end{table}
\section{Numerical Examples}
Given that some preliminary results have been given by \citet{Darrall2013,Pedgaonkar2021}, numerical examples presented in this work focus on the performance of drilling degrees of freedom offered by the consistent couple stress theory and the size dependence related topics.
\subsection{Patch Test}
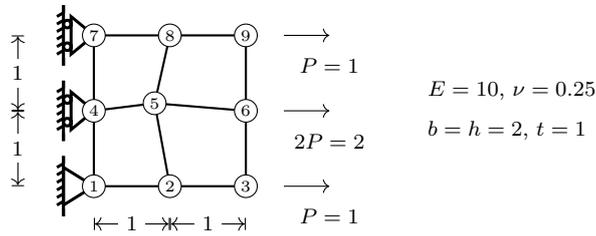
\begin{figure}[htb]
\centering\scriptsize
\begin{tikzpicture}[scale=1]
\coordinate(N1)at(0,0);
\coordinate(N2)at(1,0);
\coordinate(N3)at(2,0);
\coordinate(N4)at(0,1);
\coordinate(N5)at(.8,1.1);
\coordinate(N6)at(2,1);
\coordinate(N7)at(0,2);
\coordinate(N8)at(1,2);
\coordinate(N9)at(2,2);
\draw[thick](N1)--(N2)--(N3)--(N6)--(N9)--(N8)--(N7)--(N4)--cycle;
\draw[thick](N2)--(N5)--(N8);
\draw[thick](N4)--(N5)--(N6);
\HingeSupport[-90]{N1};
\RollerSupport[-90]{N4};
\RollerSupport[-90]{N7};
\begin{scope}[every node/.style={fill=white,circle,draw,inner sep=0,minimum size=3mm}]
\node at(N1){\tiny1};
\node at(N2){\tiny2};
\node at(N3){\tiny3};
\node at(N4){\tiny4};
\node at(N5){\tiny5};
\node at(N6){\tiny6};
\node at(N7){\tiny7};
\node at(N8){\tiny8};
\node at(N9){\tiny9};
\end{scope}
\draw[->](2.5,0)--++(.6,0)node[below=2mm]{$P=1$};
\draw[->](2.5,1)--++(.6,0)node[below=2mm]{$2P=2$};
\draw[->](2.5,2)--++(.6,0)node[below=2mm]{$P=1$};
\draw[|<->|](0,-.5)--++(1,0)node[midway,fill=white]{$1$};
\draw[|<->|](1,-.5)--++(1,0)node[midway,fill=white]{$1$};
\draw[|<->|](-1,0)--++(0,1)node[midway,fill=white]{$1$};
\draw[|<->|](-1,1)--++(0,1)node[midway,fill=white]{$1$};
\node[align=left]at(5.5,1){$E=10$, $\nu=0.25$\\[2mm]$b=h=2$, $t=1$};
\end{tikzpicture}
\caption{constant strain patch test}\label{fig:csmq_patch}
\end{figure}

As a convention adopted by the FEM community, the classic constant strain patch test that defines four elements in a rectangular panel as shown in \figref{fig:csmq_patch} is firstly presented as a preliminary to convergence. The linear displacement field can be successfully computed by using CSMQ4 elements with arbitrary location of the middle node.

Since the strain field remains constant, no couple stress would be generated in this example. As the result, whether nodal rotations are constrained has no impact on final results. Due to the same reason, any positive numbers can be chosen as the characteristic length $l$.
\subsection{Plane Ring}
To further validate the correctness of the formulation, as well as the implemented elements, the plane ring example is modelled. The corresponding analytical solution is available, as can be seen elsewhere \citep{Hadjesfandiari2011}, making it a perfect choice for the purpose of validation.

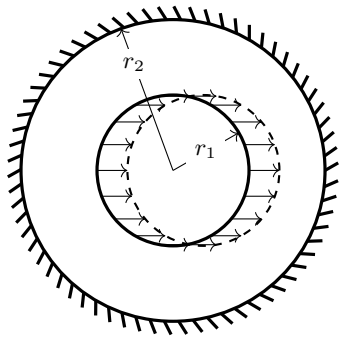
\begin{figure}[htb]
\centering\footnotesize
\begin{tikzpicture}
\draw[very thick,FFIXED](0,0)circle(2);
\draw[very thick](0,0)circle(1);
\draw[dashed,thick](0.4,0)circle(1);
\foreach \i in {0,20,...,340}{\draw[->](\i:1)--++(0.4,0);}
\draw[->](0,0)--(30:1)node[midway,fill=white]{$r_1$};
\draw[->](0,0)--(110:2)node[near end,fill=white]{$r_2$};
\end{tikzpicture}
\caption{plane ring subjected to uniform horizontal displacement}\label{fig:ring}
\end{figure}

\begin{figure}[htb]
\centering\footnotesize
\includegraphics[width=7cm]{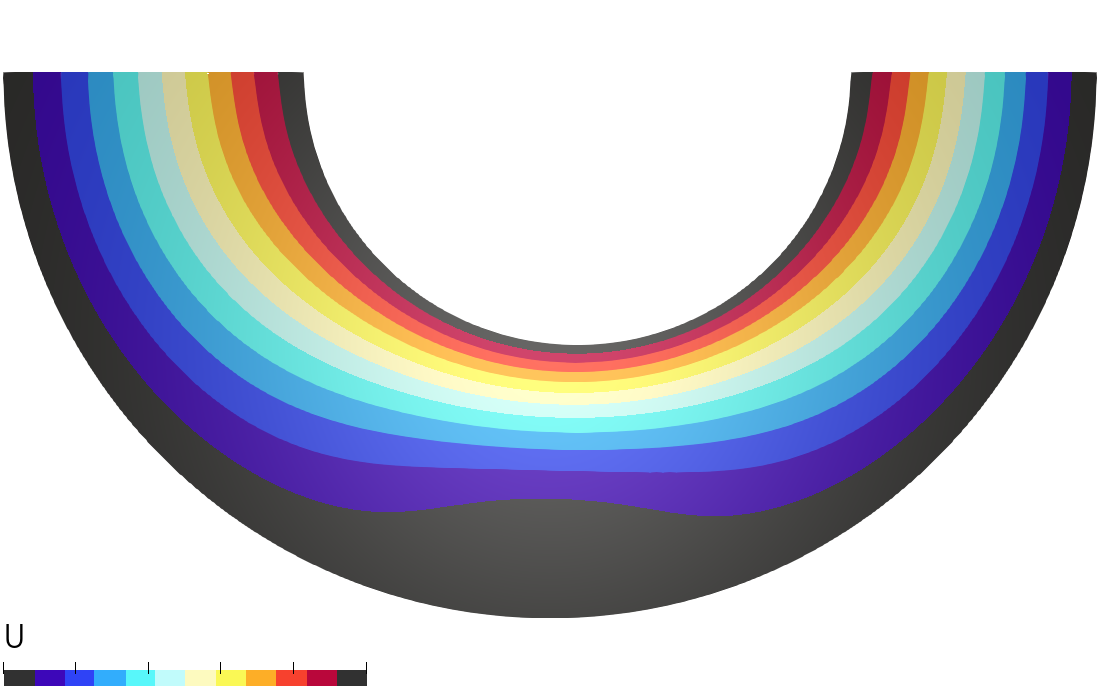}
\caption{deformed half ring model with $l=\num{0.1}$}\label{fig:deformed_ring}
\end{figure}

The ring shown in \figref{fig:ring} is subjected to plane strain condition with inner radius $r_1=1$, outer radius $r_2=2$, shear modulus $\mu=1$ and Poisson's ratio $\nu=0.4$.
The outer boundary is fixed thus translation is zero. A uniform horizontal displacement is applied to the inner boundary while the vertical displacement is constrained. As free couple traction is assumed on both boundaries, the corresponding rotation is not necessarily zero. Due to symmetry, the finite element model defines the geometry of half of the ring with a structured grid of size \num{25x100}. Three elements are tested: CSMT3, CSMQ4 and CSMQ8. The deformed model can be seen in \figref{fig:deformed_ring}.

\begin{figure}[htb]
\centering\footnotesize
\begin{subfigure}[b]{.49\textwidth}\centering
\begin{tikzpicture}[gnuplot]
\path (0.000,0.000) rectangle (6.000,4.000);
\gpcolor{rgb color={0.933,0.933,0.933}}
\gpsetlinetype{gp lt border}
\gpsetdashtype{gp dt solid}
\gpsetlinewidth{1.00}
\draw[gp path] (0.000,0.190)--(5.999,0.190);
\gpcolor{color=gp lt color border}
\draw[gp path] (0.000,0.190)--(0.180,0.190);
\draw[gp path] (5.999,0.190)--(5.819,0.190);
\node[gp node right] at (-0.184,0.190) {$0$};
\gpcolor{rgb color={0.933,0.933,0.933}}
\draw[gp path] (0.000,0.952)--(5.999,0.952);
\gpcolor{color=gp lt color border}
\draw[gp path] (0.000,0.952)--(0.180,0.952);
\draw[gp path] (5.999,0.952)--(5.819,0.952);
\node[gp node right] at (-0.184,0.952) {$0.2$};
\gpcolor{rgb color={0.933,0.933,0.933}}
\draw[gp path] (0.000,1.714)--(5.999,1.714);
\gpcolor{color=gp lt color border}
\draw[gp path] (0.000,1.714)--(0.180,1.714);
\draw[gp path] (5.999,1.714)--(5.819,1.714);
\node[gp node right] at (-0.184,1.714) {$0.4$};
\gpcolor{rgb color={0.933,0.933,0.933}}
\draw[gp path] (0.000,2.476)--(3.311,2.476);
\draw[gp path] (5.699,2.476)--(5.999,2.476);
\gpcolor{color=gp lt color border}
\draw[gp path] (0.000,2.476)--(0.180,2.476);
\draw[gp path] (5.999,2.476)--(5.819,2.476);
\node[gp node right] at (-0.184,2.476) {$0.6$};
\gpcolor{rgb color={0.933,0.933,0.933}}
\draw[gp path] (0.000,3.237)--(3.311,3.237);
\draw[gp path] (5.699,3.237)--(5.999,3.237);
\gpcolor{color=gp lt color border}
\draw[gp path] (0.000,3.237)--(0.180,3.237);
\draw[gp path] (5.999,3.237)--(5.819,3.237);
\node[gp node right] at (-0.184,3.237) {$0.8$};
\gpcolor{rgb color={0.933,0.933,0.933}}
\draw[gp path] (0.000,3.999)--(5.999,3.999);
\gpcolor{color=gp lt color border}
\draw[gp path] (0.000,3.999)--(0.180,3.999);
\draw[gp path] (5.999,3.999)--(5.819,3.999);
\node[gp node right] at (-0.184,3.999) {$1$};
\gpcolor{rgb color={0.933,0.933,0.933}}
\draw[gp path] (0.000,0.000)--(0.000,3.999);
\gpcolor{color=gp lt color border}
\draw[gp path] (0.000,0.000)--(0.000,0.180);
\draw[gp path] (0.000,3.999)--(0.000,3.819);
\node[gp node center] at (0.000,-0.308) {$1$};
\gpcolor{rgb color={0.933,0.933,0.933}}
\draw[gp path] (1.200,0.000)--(1.200,3.999);
\gpcolor{color=gp lt color border}
\draw[gp path] (1.200,0.000)--(1.200,0.180);
\draw[gp path] (1.200,3.999)--(1.200,3.819);
\node[gp node center] at (1.200,-0.308) {$1.2$};
\gpcolor{rgb color={0.933,0.933,0.933}}
\draw[gp path] (2.400,0.000)--(2.400,3.999);
\gpcolor{color=gp lt color border}
\draw[gp path] (2.400,0.000)--(2.400,0.180);
\draw[gp path] (2.400,3.999)--(2.400,3.819);
\node[gp node center] at (2.400,-0.308) {$1.4$};
\gpcolor{rgb color={0.933,0.933,0.933}}
\draw[gp path] (3.599,0.000)--(3.599,2.323);
\draw[gp path] (3.599,3.799)--(3.599,3.999);
\gpcolor{color=gp lt color border}
\draw[gp path] (3.599,0.000)--(3.599,0.180);
\draw[gp path] (3.599,3.999)--(3.599,3.819);
\node[gp node center] at (3.599,-0.308) {$1.6$};
\gpcolor{rgb color={0.933,0.933,0.933}}
\draw[gp path] (4.799,0.000)--(4.799,2.323);
\draw[gp path] (4.799,3.799)--(4.799,3.999);
\gpcolor{color=gp lt color border}
\draw[gp path] (4.799,0.000)--(4.799,0.180);
\draw[gp path] (4.799,3.999)--(4.799,3.819);
\node[gp node center] at (4.799,-0.308) {$1.8$};
\gpcolor{rgb color={0.933,0.933,0.933}}
\draw[gp path] (5.999,0.000)--(5.999,3.999);
\gpcolor{color=gp lt color border}
\draw[gp path] (5.999,0.000)--(5.999,0.180);
\draw[gp path] (5.999,3.999)--(5.999,3.819);
\node[gp node center] at (5.999,-0.308) {$2$};
\draw[gp path] (0.000,3.999)--(0.000,0.000)--(5.999,0.000)--(5.999,3.999)--cycle;
\node[gp node center,rotate=-270] at (-1.028,1.999) {normalised displacement};
\node[gp node center] at (2.999,-0.769) {radius};
\node[gp node left] at (3.311,3.614) {analytical};
\gpcolor{rgb color={0.000,0.000,0.000}}
\gpsetlinewidth{6.00}
\draw[gp path] (5.151,3.614)--(5.515,3.614);
\draw[gp path] (0.000,3.999)--(0.120,3.797)--(0.240,3.599)--(0.360,3.409)--(0.480,3.218)%
  --(0.600,3.039)--(0.720,2.864)--(0.840,2.693)--(0.960,2.533)--(1.080,2.377)--(1.200,2.224)%
  --(1.320,2.083)--(1.440,1.946)--(1.560,1.813)--(1.680,1.691)--(1.800,1.573)--(1.920,1.459)%
  --(2.040,1.352)--(2.160,1.249)--(2.280,1.154)--(2.400,1.063)--(2.520,0.979)--(2.640,0.899)%
  --(2.760,0.823)--(2.880,0.750)--(3.000,0.686)--(3.119,0.625)--(3.239,0.566)--(3.359,0.513)%
  --(3.479,0.465)--(3.599,0.420)--(3.719,0.379)--(3.839,0.342)--(3.959,0.308)--(4.079,0.279)%
  --(4.199,0.253)--(4.319,0.230)--(4.439,0.211)--(4.559,0.195)--(4.679,0.182)--(4.799,0.172)%
  --(4.919,0.165)--(5.039,0.161)--(5.159,0.159)--(5.279,0.159)--(5.399,0.161)--(5.519,0.165)%
  --(5.639,0.170)--(5.759,0.177)--(5.879,0.184)--(5.999,0.190);
\gpcolor{color=gp lt color border}
\node[gp node left] at (3.311,3.245) {CSMQ8};
\gpcolor{rgb color={0.894,0.102,0.110}}
\gpsetlinewidth{3.00}
\gpsetpointsize{7.20}
\gp3point{gp mark 6}{}{(0.000,3.999)}
\gp3point{gp mark 6}{}{(0.240,3.599)}
\gp3point{gp mark 6}{}{(0.480,3.218)}
\gp3point{gp mark 6}{}{(0.720,2.864)}
\gp3point{gp mark 6}{}{(0.960,2.533)}
\gp3point{gp mark 6}{}{(1.200,2.224)}
\gp3point{gp mark 6}{}{(1.440,1.946)}
\gp3point{gp mark 6}{}{(1.680,1.691)}
\gp3point{gp mark 6}{}{(1.920,1.459)}
\gp3point{gp mark 6}{}{(2.160,1.249)}
\gp3point{gp mark 6}{}{(2.400,1.063)}
\gp3point{gp mark 6}{}{(2.640,0.899)}
\gp3point{gp mark 6}{}{(2.880,0.750)}
\gp3point{gp mark 6}{}{(3.119,0.625)}
\gp3point{gp mark 6}{}{(3.359,0.513)}
\gp3point{gp mark 6}{}{(3.599,0.420)}
\gp3point{gp mark 6}{}{(3.839,0.342)}
\gp3point{gp mark 6}{}{(4.079,0.279)}
\gp3point{gp mark 6}{}{(4.319,0.230)}
\gp3point{gp mark 6}{}{(4.559,0.195)}
\gp3point{gp mark 6}{}{(4.799,0.172)}
\gp3point{gp mark 6}{}{(5.039,0.161)}
\gp3point{gp mark 6}{}{(5.279,0.159)}
\gp3point{gp mark 6}{}{(5.519,0.165)}
\gp3point{gp mark 6}{}{(5.759,0.176)}
\gp3point{gp mark 6}{}{(5.999,0.190)}
\gp3point{gp mark 6}{}{(5.333,3.245)}
\gpcolor{color=gp lt color border}
\node[gp node left] at (3.311,2.876) {CSMQ4};
\gpcolor{rgb color={1.000,0.498,0.000}}
\gp3point{gp mark 3}{}{(0.000,3.999)}
\gp3point{gp mark 3}{}{(0.300,3.550)}
\gp3point{gp mark 3}{}{(0.600,3.115)}
\gp3point{gp mark 3}{}{(0.900,2.708)}
\gp3point{gp mark 3}{}{(1.200,2.331)}
\gp3point{gp mark 3}{}{(1.500,1.988)}
\gp3point{gp mark 3}{}{(1.800,1.676)}
\gp3point{gp mark 3}{}{(2.100,1.398)}
\gp3point{gp mark 3}{}{(2.400,1.154)}
\gp3point{gp mark 3}{}{(2.700,0.941)}
\gp3point{gp mark 3}{}{(3.000,0.758)}
\gp3point{gp mark 3}{}{(3.299,0.606)}
\gp3point{gp mark 3}{}{(3.599,0.477)}
\gp3point{gp mark 3}{}{(3.899,0.375)}
\gp3point{gp mark 3}{}{(4.199,0.296)}
\gp3point{gp mark 3}{}{(4.499,0.239)}
\gp3point{gp mark 3}{}{(4.799,0.202)}
\gp3point{gp mark 3}{}{(5.099,0.182)}
\gp3point{gp mark 3}{}{(5.399,0.176)}
\gp3point{gp mark 3}{}{(5.699,0.181)}
\gp3point{gp mark 3}{}{(5.999,0.190)}
\gp3point{gp mark 3}{}{(5.333,2.876)}
\gpcolor{color=gp lt color border}
\node[gp node left] at (3.311,2.507) {CSMT3};
\gpcolor{rgb color={0.200,0.627,0.173}}
\gp3point{gp mark 1}{}{(0.000,3.999)}
\gp3point{gp mark 1}{}{(0.480,3.333)}
\gp3point{gp mark 1}{}{(0.960,2.700)}
\gp3point{gp mark 1}{}{(1.440,2.129)}
\gp3point{gp mark 1}{}{(1.920,1.638)}
\gp3point{gp mark 1}{}{(2.400,1.230)}
\gp3point{gp mark 1}{}{(2.880,0.895)}
\gp3point{gp mark 1}{}{(3.359,0.636)}
\gp3point{gp mark 1}{}{(3.839,0.442)}
\gp3point{gp mark 1}{}{(4.319,0.309)}
\gp3point{gp mark 1}{}{(4.799,0.228)}
\gp3point{gp mark 1}{}{(5.279,0.191)}
\gp3point{gp mark 1}{}{(5.759,0.186)}
\gp3point{gp mark 1}{}{(5.333,2.507)}
\gpcolor{color=gp lt color border}
\gpsetlinewidth{1.00}
\draw[gp path] (0.000,3.999)--(0.000,0.000)--(5.999,0.000)--(5.999,3.999)--cycle;
\gpdefrectangularnode{gp plot 1}{\pgfpoint{0.000cm}{0.000cm}}{\pgfpoint{5.999cm}{3.999cm}}
\end{tikzpicture}
\caption{$l=\num{0.1}$}
\end{subfigure}\hfill
\begin{subfigure}[b]{.49\textwidth}\centering
\begin{tikzpicture}[gnuplot]
\path (0.000,0.000) rectangle (6.000,4.000);
\gpcolor{rgb color={0.933,0.933,0.933}}
\gpsetlinetype{gp lt border}
\gpsetdashtype{gp dt solid}
\gpsetlinewidth{1.00}
\draw[gp path] (0.000,0.190)--(5.999,0.190);
\gpcolor{color=gp lt color border}
\draw[gp path] (0.000,0.190)--(0.180,0.190);
\draw[gp path] (5.999,0.190)--(5.819,0.190);
\node[gp node right] at (-0.184,0.190) {$0$};
\gpcolor{rgb color={0.933,0.933,0.933}}
\draw[gp path] (0.000,0.952)--(5.999,0.952);
\gpcolor{color=gp lt color border}
\draw[gp path] (0.000,0.952)--(0.180,0.952);
\draw[gp path] (5.999,0.952)--(5.819,0.952);
\node[gp node right] at (-0.184,0.952) {$0.2$};
\gpcolor{rgb color={0.933,0.933,0.933}}
\draw[gp path] (0.000,1.714)--(5.999,1.714);
\gpcolor{color=gp lt color border}
\draw[gp path] (0.000,1.714)--(0.180,1.714);
\draw[gp path] (5.999,1.714)--(5.819,1.714);
\node[gp node right] at (-0.184,1.714) {$0.4$};
\gpcolor{rgb color={0.933,0.933,0.933}}
\draw[gp path] (0.000,2.476)--(3.311,2.476);
\draw[gp path] (5.699,2.476)--(5.999,2.476);
\gpcolor{color=gp lt color border}
\draw[gp path] (0.000,2.476)--(0.180,2.476);
\draw[gp path] (5.999,2.476)--(5.819,2.476);
\node[gp node right] at (-0.184,2.476) {$0.6$};
\gpcolor{rgb color={0.933,0.933,0.933}}
\draw[gp path] (0.000,3.237)--(3.311,3.237);
\draw[gp path] (5.699,3.237)--(5.999,3.237);
\gpcolor{color=gp lt color border}
\draw[gp path] (0.000,3.237)--(0.180,3.237);
\draw[gp path] (5.999,3.237)--(5.819,3.237);
\node[gp node right] at (-0.184,3.237) {$0.8$};
\gpcolor{rgb color={0.933,0.933,0.933}}
\draw[gp path] (0.000,3.999)--(5.999,3.999);
\gpcolor{color=gp lt color border}
\draw[gp path] (0.000,3.999)--(0.180,3.999);
\draw[gp path] (5.999,3.999)--(5.819,3.999);
\node[gp node right] at (-0.184,3.999) {$1$};
\gpcolor{rgb color={0.933,0.933,0.933}}
\draw[gp path] (0.000,0.000)--(0.000,3.999);
\gpcolor{color=gp lt color border}
\draw[gp path] (0.000,0.000)--(0.000,0.180);
\draw[gp path] (0.000,3.999)--(0.000,3.819);
\node[gp node center] at (0.000,-0.308) {$1$};
\gpcolor{rgb color={0.933,0.933,0.933}}
\draw[gp path] (1.200,0.000)--(1.200,3.999);
\gpcolor{color=gp lt color border}
\draw[gp path] (1.200,0.000)--(1.200,0.180);
\draw[gp path] (1.200,3.999)--(1.200,3.819);
\node[gp node center] at (1.200,-0.308) {$1.2$};
\gpcolor{rgb color={0.933,0.933,0.933}}
\draw[gp path] (2.400,0.000)--(2.400,3.999);
\gpcolor{color=gp lt color border}
\draw[gp path] (2.400,0.000)--(2.400,0.180);
\draw[gp path] (2.400,3.999)--(2.400,3.819);
\node[gp node center] at (2.400,-0.308) {$1.4$};
\gpcolor{rgb color={0.933,0.933,0.933}}
\draw[gp path] (3.599,0.000)--(3.599,2.323);
\draw[gp path] (3.599,3.799)--(3.599,3.999);
\gpcolor{color=gp lt color border}
\draw[gp path] (3.599,0.000)--(3.599,0.180);
\draw[gp path] (3.599,3.999)--(3.599,3.819);
\node[gp node center] at (3.599,-0.308) {$1.6$};
\gpcolor{rgb color={0.933,0.933,0.933}}
\draw[gp path] (4.799,0.000)--(4.799,2.323);
\draw[gp path] (4.799,3.799)--(4.799,3.999);
\gpcolor{color=gp lt color border}
\draw[gp path] (4.799,0.000)--(4.799,0.180);
\draw[gp path] (4.799,3.999)--(4.799,3.819);
\node[gp node center] at (4.799,-0.308) {$1.8$};
\gpcolor{rgb color={0.933,0.933,0.933}}
\draw[gp path] (5.999,0.000)--(5.999,3.999);
\gpcolor{color=gp lt color border}
\draw[gp path] (5.999,0.000)--(5.999,0.180);
\draw[gp path] (5.999,3.999)--(5.999,3.819);
\node[gp node center] at (5.999,-0.308) {$2$};
\draw[gp path] (0.000,3.999)--(0.000,0.000)--(5.999,0.000)--(5.999,3.999)--cycle;
\node[gp node center,rotate=-270] at (-1.028,1.999) {normalised displacement};
\node[gp node center] at (2.999,-0.769) {radius};
\node[gp node left] at (3.311,3.614) {analytical};
\gpcolor{rgb color={0.000,0.000,0.000}}
\gpsetlinewidth{6.00}
\draw[gp path] (5.151,3.614)--(5.515,3.614);
\draw[gp path] (0.000,3.999)--(0.120,3.885)--(0.240,3.774)--(0.360,3.664)--(0.480,3.553)%
  --(0.600,3.447)--(0.720,3.340)--(0.840,3.233)--(0.960,3.127)--(1.080,3.024)--(1.200,2.925)%
  --(1.320,2.826)--(1.440,2.727)--(1.560,2.632)--(1.680,2.537)--(1.800,2.445)--(1.920,2.354)%
  --(2.040,2.266)--(2.160,2.179)--(2.280,2.095)--(2.400,2.011)--(2.520,1.931)--(2.640,1.851)%
  --(2.760,1.775)--(2.880,1.699)--(3.000,1.622)--(3.119,1.550)--(3.239,1.482)--(3.359,1.413)%
  --(3.479,1.344)--(3.599,1.280)--(3.719,1.215)--(3.839,1.150)--(3.959,1.089)--(4.079,1.028)%
  --(4.199,0.971)--(4.319,0.914)--(4.439,0.857)--(4.559,0.804)--(4.679,0.746)--(4.799,0.693)%
  --(4.919,0.644)--(5.039,0.590)--(5.159,0.539)--(5.279,0.488)--(5.399,0.438)--(5.519,0.388)%
  --(5.639,0.339)--(5.759,0.289)--(5.879,0.240)--(5.999,0.190);
\gpcolor{color=gp lt color border}
\node[gp node left] at (3.311,3.245) {CSMQ8};
\gpcolor{rgb color={0.894,0.102,0.110}}
\gpsetlinewidth{3.00}
\gpsetpointsize{7.20}
\gp3point{gp mark 6}{}{(0.000,3.999)}
\gp3point{gp mark 6}{}{(0.240,3.774)}
\gp3point{gp mark 6}{}{(0.480,3.553)}
\gp3point{gp mark 6}{}{(0.720,3.336)}
\gp3point{gp mark 6}{}{(0.960,3.127)}
\gp3point{gp mark 6}{}{(1.200,2.921)}
\gp3point{gp mark 6}{}{(1.440,2.727)}
\gp3point{gp mark 6}{}{(1.680,2.537)}
\gp3point{gp mark 6}{}{(1.920,2.354)}
\gp3point{gp mark 6}{}{(2.160,2.179)}
\gp3point{gp mark 6}{}{(2.400,2.011)}
\gp3point{gp mark 6}{}{(2.640,1.847)}
\gp3point{gp mark 6}{}{(2.880,1.695)}
\gp3point{gp mark 6}{}{(3.119,1.550)}
\gp3point{gp mark 6}{}{(3.359,1.409)}
\gp3point{gp mark 6}{}{(3.599,1.276)}
\gp3point{gp mark 6}{}{(3.839,1.150)}
\gp3point{gp mark 6}{}{(4.079,1.028)}
\gp3point{gp mark 6}{}{(4.319,0.914)}
\gp3point{gp mark 6}{}{(4.559,0.800)}
\gp3point{gp mark 6}{}{(4.799,0.693)}
\gp3point{gp mark 6}{}{(5.039,0.590)}
\gp3point{gp mark 6}{}{(5.279,0.487)}
\gp3point{gp mark 6}{}{(5.519,0.387)}
\gp3point{gp mark 6}{}{(5.759,0.289)}
\gp3point{gp mark 6}{}{(5.999,0.190)}
\gp3point{gp mark 6}{}{(5.333,3.245)}
\gpcolor{color=gp lt color border}
\node[gp node left] at (3.311,2.876) {CSMQ4};
\gpcolor{rgb color={1.000,0.498,0.000}}
\gp3point{gp mark 3}{}{(0.000,3.999)}
\gp3point{gp mark 3}{}{(0.300,3.759)}
\gp3point{gp mark 3}{}{(0.600,3.519)}
\gp3point{gp mark 3}{}{(0.900,3.283)}
\gp3point{gp mark 3}{}{(1.200,3.051)}
\gp3point{gp mark 3}{}{(1.500,2.822)}
\gp3point{gp mark 3}{}{(1.800,2.605)}
\gp3point{gp mark 3}{}{(2.100,2.388)}
\gp3point{gp mark 3}{}{(2.400,2.182)}
\gp3point{gp mark 3}{}{(2.700,1.984)}
\gp3point{gp mark 3}{}{(3.000,1.790)}
\gp3point{gp mark 3}{}{(3.299,1.603)}
\gp3point{gp mark 3}{}{(3.599,1.428)}
\gp3point{gp mark 3}{}{(3.899,1.257)}
\gp3point{gp mark 3}{}{(4.199,1.089)}
\gp3point{gp mark 3}{}{(4.499,0.933)}
\gp3point{gp mark 3}{}{(4.799,0.777)}
\gp3point{gp mark 3}{}{(5.099,0.628)}
\gp3point{gp mark 3}{}{(5.399,0.479)}
\gp3point{gp mark 3}{}{(5.699,0.335)}
\gp3point{gp mark 3}{}{(5.999,0.190)}
\gp3point{gp mark 3}{}{(5.333,2.876)}
\gpcolor{color=gp lt color border}
\node[gp node left] at (3.311,2.507) {CSMT3};
\gpcolor{rgb color={0.200,0.627,0.173}}
\gp3point{gp mark 1}{}{(0.000,3.999)}
\gp3point{gp mark 1}{}{(0.480,3.630)}
\gp3point{gp mark 1}{}{(0.960,3.268)}
\gp3point{gp mark 1}{}{(1.440,2.914)}
\gp3point{gp mark 1}{}{(1.920,2.571)}
\gp3point{gp mark 1}{}{(2.400,2.239)}
\gp3point{gp mark 1}{}{(2.880,1.927)}
\gp3point{gp mark 1}{}{(3.359,1.630)}
\gp3point{gp mark 1}{}{(3.839,1.344)}
\gp3point{gp mark 1}{}{(4.319,1.074)}
\gp3point{gp mark 1}{}{(4.799,0.815)}
\gp3point{gp mark 1}{}{(5.279,0.564)}
\gp3point{gp mark 1}{}{(5.759,0.314)}
\gp3point{gp mark 1}{}{(5.333,2.507)}
\gpcolor{color=gp lt color border}
\gpsetlinewidth{1.00}
\draw[gp path] (0.000,3.999)--(0.000,0.000)--(5.999,0.000)--(5.999,3.999)--cycle;
\gpdefrectangularnode{gp plot 1}{\pgfpoint{0.000cm}{0.000cm}}{\pgfpoint{5.999cm}{3.999cm}}
\end{tikzpicture}
\caption{$l=\num{0.5}$}
\end{subfigure}
\caption{$u_\theta$ along the vertical center line}\label{fig:u_theta}
\end{figure}
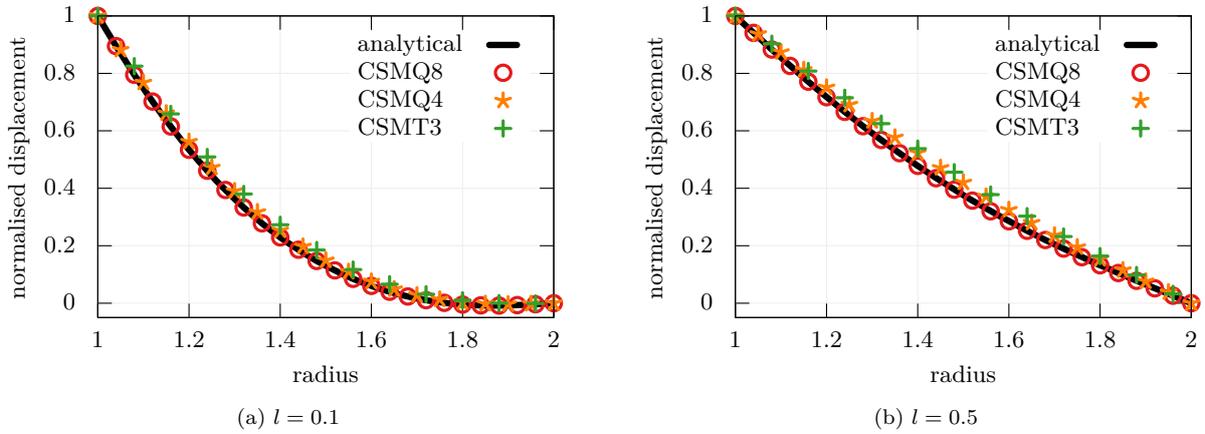

Numerical results of the transverse displacement $u_\theta$ along the vertical center line are presented in \figref{fig:u_theta} with analytical solution obtained by using proper $l$.
The numerical solution matches the analytical one well, indicating the formulation/implementation is correct. Since no dedicated optimisation is designed to improve the performance apart from the mixed formulation, relatively dense mesh grids are required to reduce error. The second order quadrilateral CSMQ8 shows higher accuracy as can be inferred based on the experience with the Cauchy framework.
\subsection{Membrane--Beam Joint}
From the previous two examples, it can be concluded that the formulation developed in this work and the corresponding implementation are both correct. Since the drilling degree of freedom defined in the consistent couple stress theory has an energetic conjugate --- curl of couple stress $\nabla\times\mb{\mu}$ as can be seen in \eqsref{eq:functional2}, size effect can be accounted for by its nature.

It is observed in the previous work \citep{Chang2020} that drilling DoF formulated within the Cauchy framework tends to be mesh sensitive --- in-plane rotational stiffness decreases with increasing mesh size. This leads to unreliable predictions of in-plane rotation and its force counterpart. This is mainly caused by the fact that within the Cauchy framework, the drilling displacement lacks the corresponding energetic conjugate and has to be converted to translational displacement via some technique. Fundamentally, the energy potential is still governed by conventional displacement/strain and force/stress. However, the couple stress theory may offer an alternative to define mesh objective drilling DoFs.

In this example, the mesh dependency of drilling DoFs is investigated. Although it is possible to directly apply loads to drilling DoFs as `moments', a membrane--beam joint, which is frequently encountered in structural engineering, is used instead.

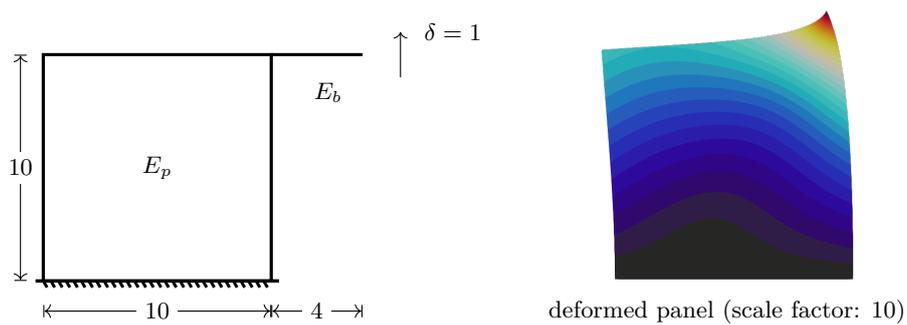
\begin{figure}[htb]
\centering\footnotesize
\begin{tikzpicture}
\FixedSupport{1.5,0}{4}
\draw[very thick](0,0)rectangle(3,3);
\draw[very thick](3,3)--(4.2,3);
\draw[->](4.7,2.7)--(4.7,3.3)node[right=2mm]{$\delta=1$};
\node[]at(1.5,1.5){$E_p$};
\node[]at(3.75,2.5){$E_b$};
\draw[|<->|](-.3,0)--++(0,3)node[midway,fill=white]{$10$};
\draw[|<->|](0,-.4)--++(3,0)node[midway,fill=white]{$10$};
\draw[|<->|](3,-.4)--++(1.2,0)node[midway,fill=white]{$4$};
\node[]at(9,1.8){\includegraphics[height=4cm]{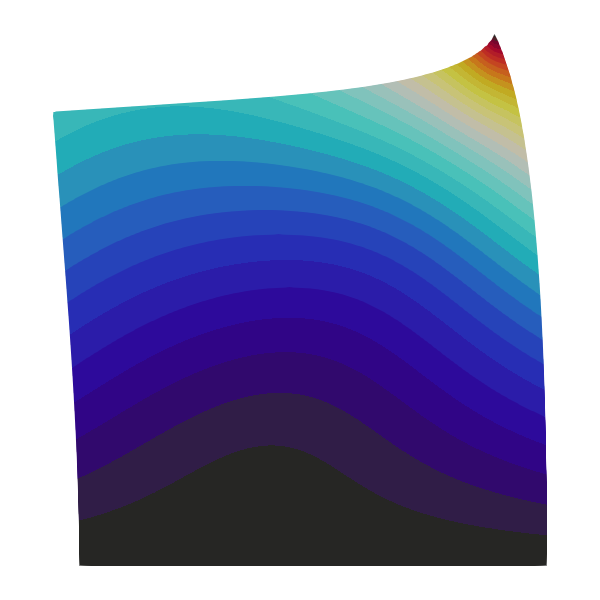}};
\node[]at(9,-.4){deformed panel (scale factor: \num{10})};
\end{tikzpicture}
\caption{panel joint with attached beam}\label{fig:joint_example}
\end{figure}

The model shown in \figref{fig:joint_example} is subjected to plane stress condition and consists of a square panel of size \num{10x10} and a beam of length \num{4}. The height of the beam section is \num{1} while the thickness along the $z$-axis (perpendicular to the plane of paper) is set to unity.
A unit displacement is applied to the free end of beam and the corresponding resistance is recorded for comparison. The panel is modelled by a structured grid with various numbers of different membrane elements (Simplified Generalised Conforming Mixed Quadrilateral with Gauss quadrature, SGCMQG \citep{Chang2020} and CSMQ4) along each direction while the beam is modelled by a single elastic Euler--Bernoulli beam element. Here, SGCMQG element is a recently proposed high-performing element with drilling DoFs constructed within the Cauchy framework. It is used as a reference to represent the typical behaviour of Cauchy--type drilling DoFs, which exhibits significant mesh dependency as can be seen later. With such a setup, a mixed action (rotation and translation) is experienced by the panel.

\begin{figure}[ht]
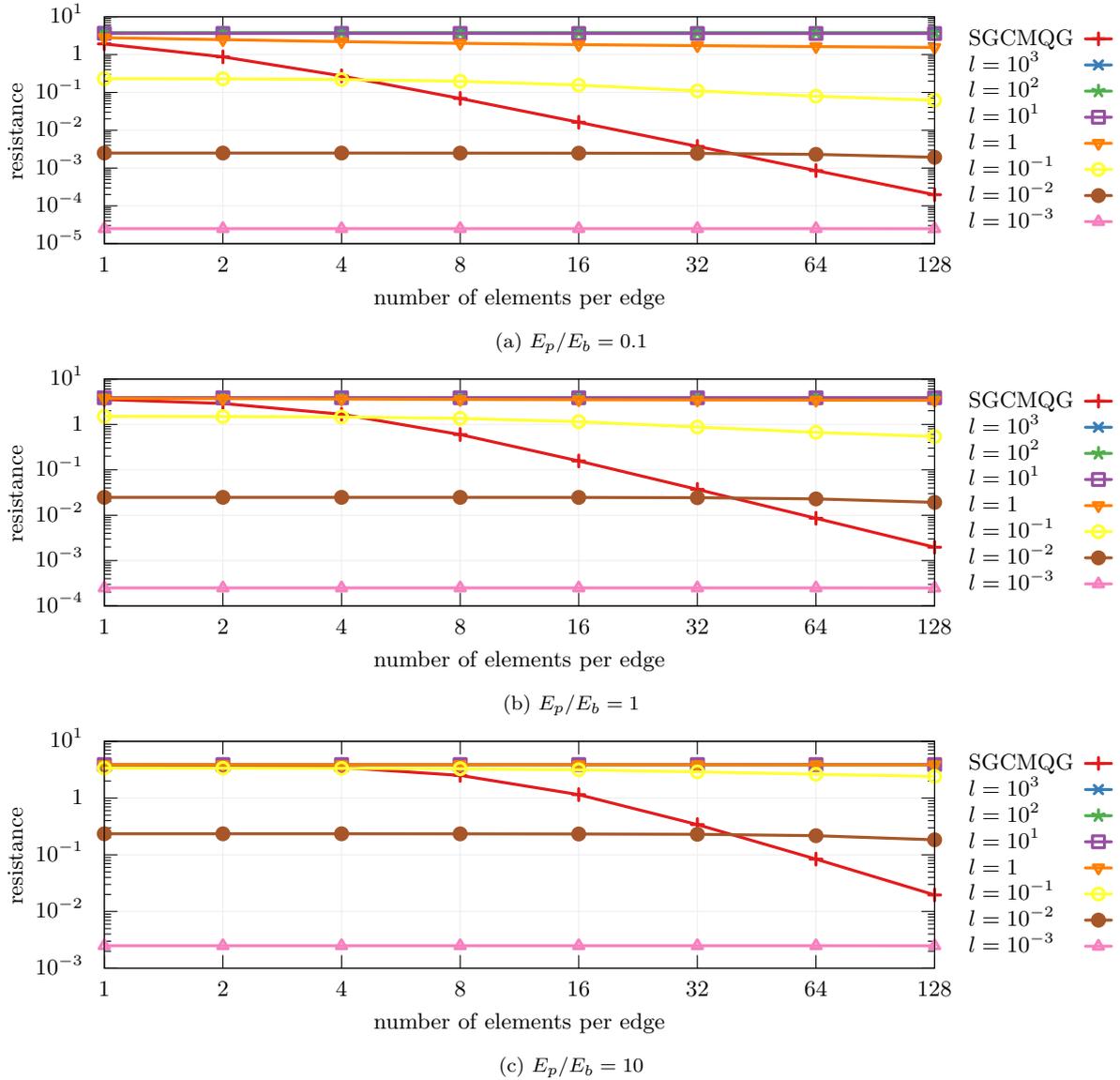

\centering\footnotesize
\begin{subfigure}[b]{.99\textwidth}\centering
\input{PIC/JOINTA}
\caption{$E_p/E_b=0.1$}
\end{subfigure}
\begin{subfigure}[b]{.99\textwidth}\centering
\input{PIC/JOINTB}
\caption{$E_p/E_b=1$}
\end{subfigure}
\begin{subfigure}[b]{.99\textwidth}\centering
\input{PIC/JOINTC}
\caption{$E_p/E_b=10$}
\end{subfigure}
\caption{beam end resistance subjected to unit vertical displacement}\label{fig:joint}
\end{figure}
Three pairs of panel and beam elastic moduli are chosen to represent different deformation patterns. The analysis is performed with characteristic length $l$ ranging from \num{E3} to \num{E-3} for CSMQ4. Numerical results are presented in \figref{fig:joint}. The following observations can be made.
\begin{enumerate}
\item When either $l$ or $E_p/E_b$ is large, the panel is sufficiently rigid so that deformation is mainly contributed by the beam, resulting in a cantilever like structure, the corresponding resistance is close to the upper bound \num{3.9063} computed via a cantilever model with the Euler--Bernoulli theory.
\item With all three different $E_p/E_b$ ratios, for the Cauchy theory based SGCMQG element, the corresponding rotational stiffness decreases (thus deformation localises and resistance drops) with refined meshes. On the contrary, CSMQ4 is less sensitive to mesh size. With fixed $l$ values, the resistance given by CSMQ4 tends to be at a constant level on a logarithmic scale and its variation is significantly smaller than that of SGCMQG.
\item With CSMQ4, mesh refinement has the most noticeable impact on numerical results when $l=\numrange{0.1}{1}$. From \figref{fig:joint}, the resistance drops with denser mesh grids for $l=0.1$ in all three cases.
\end{enumerate}

It can be further inferred that with sufficient number of elements, both SGCMQG and CSMQ4 (with a sufficiently small $l$) would converge to the same solution in which deformation is localised in a small region around the top right corner.

Given that SGCMQG element can converge to analytical solutions as can be seen in the previous example, it is no doubt that the consistent couple stress theory provides an appealing (and maybe more realistic) approach when it comes to similar problems involving in-plane rotations.

However, the determination of characteristic length $l$ could be a challenging task. The corresponding experimental theories shall be investigated in the future. With a properly defined $l$, it is believed that the implemented couple stress elements can model in-plane rotation with acceptable accuracy. As discussed by others \citep[e.g.,][]{Darrall2013}, a large $l$ would suppress bending action, resulting in a shear dominated deformation pattern.
\subsection{Inelastic Response}
From the last example, it can be seen that the consistent couple stress theory is able to produce mesh independent, or at least less mesh sensitive, results with linear elastic materials. However, the mesh/size dependency issue may also arise due to other reasons, such as inelastic material behaviour. Noting that in \eqsref{eq:stored_energy} the stored energy function can be split into two parts and each is independent from the other, it would be interesting to investigate the performance of finite elements for inelastic softening response with the assist of couple stress.

\begin{figure}[htb]
\centering\footnotesize
\begin{tikzpicture}
\def\a{2.0};\def\b{\a}
\draw[very thick](-\b,-\a)rectangle(\b,\a)(0,0)circle(\a/4);
\draw[|<->|](-\b,-1.3*\a)--++(2*\b,0)node[midway,fill=white]{$H=\num{10}$};
\draw[|<->|](-1.5*\b,-\a)--++(0,2*\a)node[midway,fill=white,rotate=90]{$L=\num{10}$};
\foreach\x in{-2,-1.5,-1,-.5,0,.5,1,1.5,2}{
\draw[->](\b,\x)--++(.5,0);
}
\HingeSupport[-90]{-\b,-\a};
\RollerSupport[-90]{-\b,-.5*\a};
\RollerSupport[-90]{-\b,0};
\RollerSupport[-90]{-\b,.5*\a};
\RollerSupport[-90]{-\b,\a};
\node at(1.4*\b,0){$u$};
\draw[|<->|](-\a/1.8,-\a/4)--++(0,\a/2)node[midway,fill=white,inner sep=0]{$2a=2$};
\node[align=center]at(.5,-1.3){$E=\num{1000}$ $\nu=\num{0.2}$\\$\sigma_y=\num{1}$ $b=\num{-0.02}$\\$t=\num{1}$};
\CoorOrigin{0,0}
\end{tikzpicture}
\caption{square plate with circular hole}\label{fig:plate_with_hole}
\end{figure}
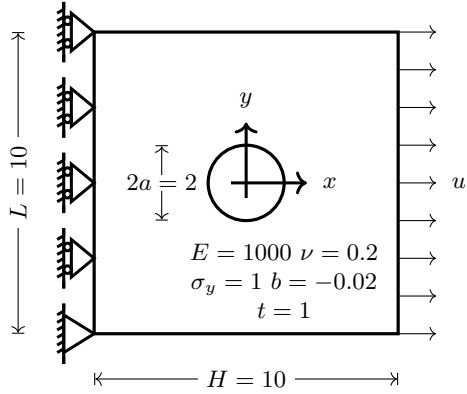
In this example, it is assumed that the couple stress part of $W$ remains `elastic' so that $\eta$ is a constant. Such an assumption is, at least, reasonable for some polymers, interested readers can refer to, for example, the work by \citet{Alisafaei2016} and the references therein. The uniform extension of a square plate with a hole of unit radius, depicted in \figref{fig:plate_with_hole}, subjected to plane stress condition is simulated with a linear softening isotropic J2 plasticity model. The following material parameters are chosen: elastic modulus $E=\num{1000}$, Poisson's ratio $\nu=\num{0.2}$, yield stress $\sigma_y=\num{1}$, hardening ratio $b=\num{-0.02}$. The thickness $t$ is set to unity.
%
%

\begin{figure}[ht]
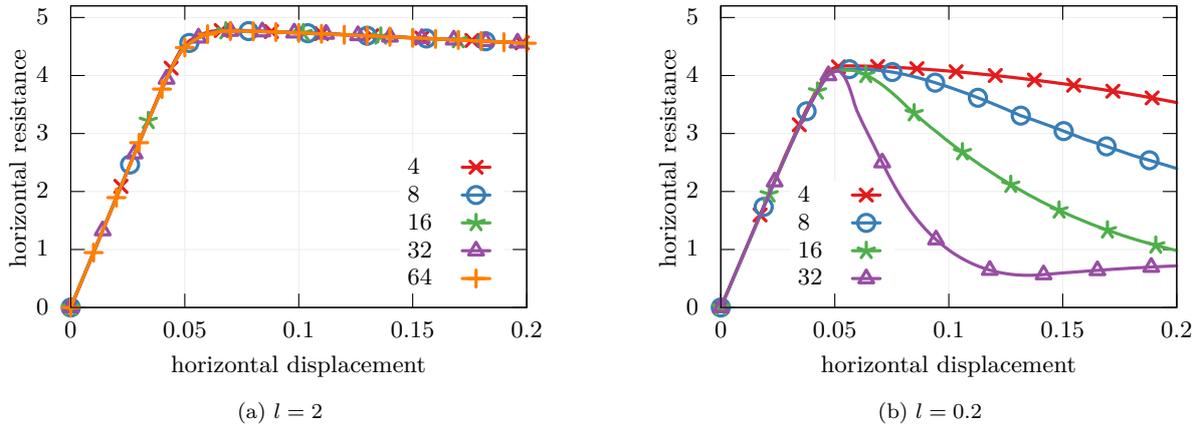

\centering\footnotesize
\begin{subfigure}[b]{.48\textwidth}\centering
\input{PIC/PLATEA}
\caption{$l=\num{2}$}\label{fig:plate_resistance_large_l}
\end{subfigure}\hfill
\begin{subfigure}[b]{.48\textwidth}\centering
\input{PIC/PLATEB}
\caption{$l=\num{0.2}$}\label{fig:plate_resistance_small_l}
\end{subfigure}
\caption{horizontal resistance of plate with various numbers of elements per edge and different characteristic lengths}
\end{figure}

From the last example, one could infer that the plate would exhibit a uniform tension like response with sufficiently large $l$ values, as in which case the last term in \eqsref{eq:stored_energy} serves as a penalty term that suppresses bending action. As can be seen in \figref{fig:plate_resistance_large_l}, with $l=\num{2}$, the displacement--resistance response of the plate resembles the strain--stress constitutive relationship regardless of element density.
However, mesh dependency cannot be resolved by the consistent couple stress theory, at least under the assumption that the couple part is `elastic'. For moderate and small $l$ values shown in \figref{fig:plate_resistance_small_l}, different mesh sizes lead to different softening curves.

\begin{figure}[ht]
\centering\footnotesize
\begin{subfigure}[b]{.33\textwidth}\centering
\includegraphics[width=4cm]{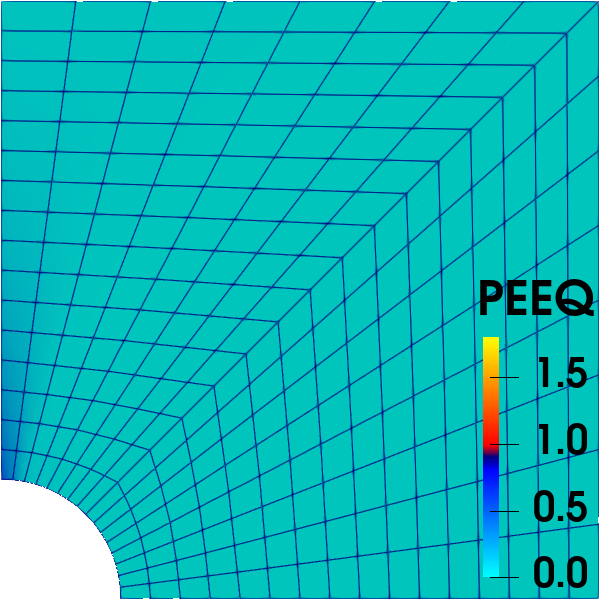}
\caption{$u=\num{0.08}$}
\end{subfigure}\hfill
\begin{subfigure}[b]{.33\textwidth}\centering
\includegraphics[width=4cm]{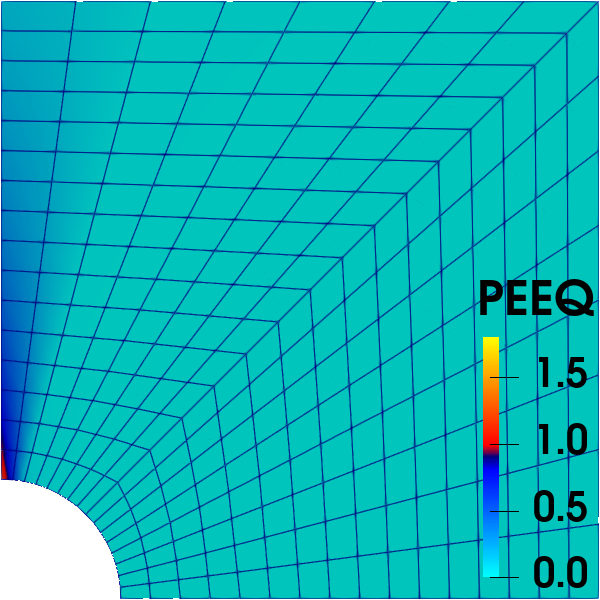}
\caption{$u=\num{0.14}$}
\end{subfigure}\hfill
\begin{subfigure}[b]{.33\textwidth}\centering
\includegraphics[width=4cm]{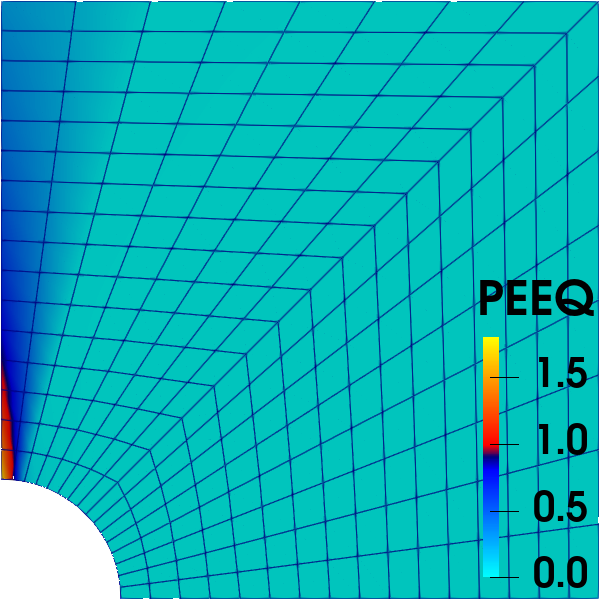}
\caption{$u=\num{0.20}$}
\end{subfigure}
\caption{distribution of equivalent plastic strain (16 elements per edge)}\label{fig:peeq16}
\begin{subfigure}[b]{.33\textwidth}\centering
\includegraphics[width=4cm]{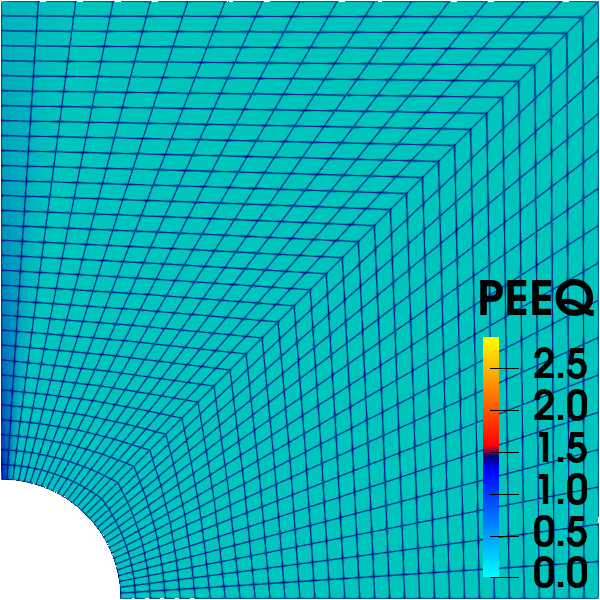}
\caption{$u=\num{0.08}$}
\end{subfigure}\hfill
\begin{subfigure}[b]{.33\textwidth}\centering
\includegraphics[width=4cm]{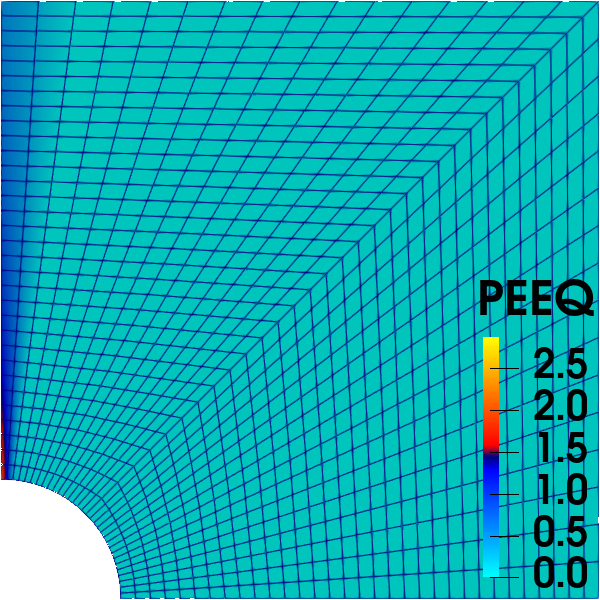}
\caption{$u=\num{0.14}$}
\end{subfigure}\hfill
\begin{subfigure}[b]{.33\textwidth}\centering
\includegraphics[width=4cm]{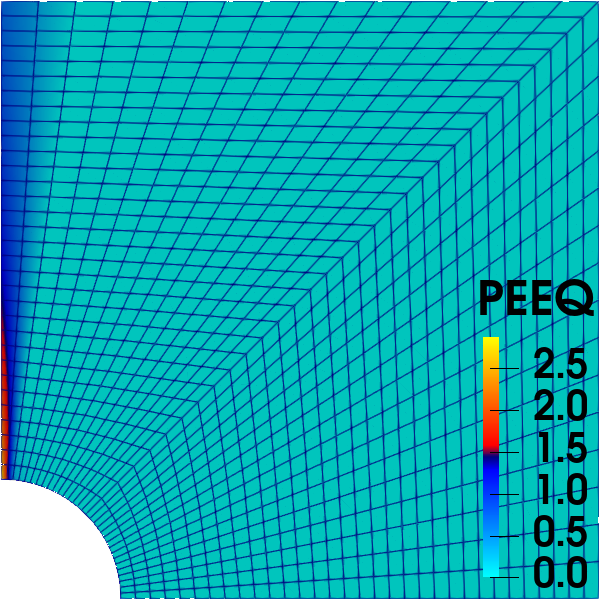}
\caption{$u=\num{0.20}$}
\end{subfigure}
\caption{distribution of equivalent plastic strain (32 elements per edge)}\label{fig:peeq32}
\end{figure}

The distribution of equivalent plastic strain (PEEQ in ABAQUS notation) is depicted in \figref{fig:peeq16} and \figref{fig:peeq32}.
Localisation of deformation, which cannot be alleviated by the introduction of couple stress field, can be identified with denser mesh grids. As the result, mesh objective results can hardly be obtained if $l$ is not large enough. Given that the consistent couple stress theory itself does not provide any regularisation of material response, such results can be predicted. A combination of couple stress theory and phase field theory \citep{Bourdin2008} can be an appealing approach in terms of modelling mesh objective softening response.

From numerical experiments, it can be inferred that $l$ seems to be more than a material constant. The geometry of the problem domain may have some effect on simulation results. Thus, the determination of $l$ may take the inherent material property (if there exists one) and the distance to the boundary of the problem domain (or some other similar characteristics) into consideration.
\section{Conclusions}
Based on the consistent couple stress theory, a mixed variational theorem with six independent fields is developed in this work. The proposed variational theorem is universal and can be converted to other existing theorems. It can also be used as a starting point to develop more advanced finite elements incorporating other techniques. A simple, efficient, local iteration free solution strategy is also proposed. With such a solution strategy, the computation cost is minimised. Examples of membrane element formulation, including three-node triangle CSMT3 and four-node quadrilateral CSMQ4, are presented. Other types of elements can be constructed following a similar procedure. Although the majority of existing constitutive models are developed based on the Cauchy theory, the proposed formulation is sufficiently versatile to accommodate both elastic and inelastic materials.

The performance of drilling degrees of freedom is evaluated via numerical examples. It is shown in the beam--panel joint example that the couple stress theory offers a better and more accurate definition of rotation (in-plane rotation in particular) which is significantly less sensitive to mesh grid size in comparison with the Cauchy theory. This increases the reliability of the derived elements and allows analysts to perform numerical analysis of models consisting of elements with various types of degrees of freedom, for example, membranes/shells plus beams. Once $l$ is properly determined, more accurate results can be obtained via the couple stress theory.

As can be seen in the inelastic example, under the assumption that the couple stress part remains elastic, mesh dependency is not significant when deformation is dominated by shear. But with decreasing $l$, deformation localises which leads to mesh dependent results. It could thus be concluded that the introduction of couple stress alone is not complete to address the mesh dependency issue in softening/fracture related problems. If $l$ is strictly defined as a material constant, let it be called material characteristic length, then there exists a mesh characteristic length that shall be used to regularise material response to achieve an objective energy release rate, in order to resolve the mesh dependency issue. Since they serve different purposes, they may not be the same quantity.

The proposed elements have been implemented in \texttt{suanPan} \citep{Chang2021}. Sample model scripts can be found online\footnote{Link to the repository would be added after acceptance}.

\bibliography{FORMULATION}
\end{document}